\documentclass[10pt,onecolumn,letterpaper]{article}
\usepackage[authoryear,sort]{natbib}
\usepackage{filecontents}
\usepackage{graphicx}
\usepackage{amssymb}
\usepackage[pdfpagelabels=true,plainpages=false,colorlinks=true,
linkcolor=blue,citecolor=blue,urlcolor=blue]{hyperref}
\usepackage[latin1]{inputenc}
\usepackage{algorithm}
\usepackage{algorithmic}
\usepackage{amsfonts}
\usepackage{amsmath}
\usepackage{cite}
\newcommand{\be}{\begin{equation}}
	\newcommand{\ee}{\end{equation}}
\newcommand{\bea}{\begin{eqnarray}}
	\newcommand{\eea}{\end{eqnarray}}
\usepackage{caption}
\usepackage{subcaption}

\usepackage{graphicx,epsfig}
\usepackage{enumerate}
\usepackage{setspace}
\usepackage{amssymb}
\usepackage{amsmath}
\usepackage{amsthm}
\usepackage{appendix}
\newtheorem{alg1}{\textbf{Algorithm}}[section]
\numberwithin{equation}{section}
\newtheorem{thm}{Theorem}[section]
\theoremstyle{definition}
\newtheorem{dfn}{Definition}[section]
\theoremstyle{Remark}
\newtheorem{note}{Note}[dfn]
\theoremstyle{lemma}
\newtheorem{rem}{Remark}[section]
\newtheorem{lem}{Lemma}[section]

\theoremstyle{example}
\newtheorem{example}{Example}[section]
\begin{document}

\title{Solution of Uncertain Multiobjective Optimization Problems by Using Nonlinear Conjugate Gradient Method}
\author{Shubham Kumar\footnote{IIITDM Jabalpur, India, E-mail: kumarshubham3394@gmail.com},
	Nihar Kumar Mahato \footnote{IIITDM Jabalpur, India, E-mail: nihar@iiitdmj.ac.in, corresponding author}
	Debdas Ghosh \footnote{IIT(BHU), India, E-mail: debdas.mat@iitbhu.ac.in},
	\\}
\date{}
\maketitle
\begin{abstract}
This paper introduces a nonlinear conjugate gradient method (NCGM) for addressing the robust counterpart of uncertain multiobjective optimization problems (UMOPs). Here, the robust counterpart is defined as the minimum across objective-wise worst-case scenarios. There are some drawbacks to using scalarization techniques to solve the robust counterparts of UMOPs, such as the pre-specification and restrictions of weights, and function importance that is unknown beforehand. NCGM  is free from any kind of priori chosen scalars
or ordering information of objective functions as accepted in scalarization methods. With the help of NCGM, we determine the critical point for the robust counterpart of UMOP, which is the robust critical point for UMOP. To tackle this robust counterpart using the NCGM, the approach involves constructing and solving a subproblem to determine a descent direction. Subsequently, a new direction is derived based on parameter selection methods such as Fletcher-Reeves, conjugate descent, Dai-Yuan, Polak-Ribi$\grave{e}$re-Polyak, and Hestenes-Stiefel. An Armijo-type inexact line search is employed to identify an appropriate step length. Utilizing descent direction and step length, a sequence is generated, and convergence of the proposed method is established. The effectiveness of the proposed method is verified and compared against an existing method using a set of test problems.
\textit{\bf 2020 Math Subject Classifications:} 90C25; 90C26; 90C29; 90C30
\end{abstract}

\section{Introduction}
Numerous factors make it difficult to apply mathematical optimization techniques to real-world problems. The multiobjective nature of the majority of real-world problems presents one significant challenge. Additionally, the input data is often unknown or subject to change, further complicating the process. As a result, finding suitable solutions for practical applications becomes a complex task.
\par In the realm of mathematical optimization, it is imperative to account for uncertainties inherent in the problem formulation when applying solutions to real-world scenarios. A critical aspect of this process involves analyzing the sensitivity of an optimal solution to variations in input data. In the literature on uncertain single-objective optimization problems (USOPs), this crucial analysis is typically conducted as a post-optimization step, commonly referred to as sensitivity analysis. For a thorough exploration of this subject, Saltelli et al. provide an insightful and comprehensive overview in \citet{saltelli2000scott}.
\par Stochastic optimization (SO) and robust optimization (RO) provide tools to directly address uncertainty during the optimization stage, in contrast to sensitivity analysis. In SO, probabilistic information about various potential realizations of uncertain input data is either provided or assumed. The goal then shifts to optimizing parameters like the mean objective value across all scenarios, the variability among scenarios, or a combination of both. For a thorough introduction and overview of the various approaches that make up stochastic optimization, see Birge and Louveaux \citet{birge2011introduction}. 
\par By not relying on any probabilistic knowledge about various scenarios of uncertain input data, RO adopts a different strategy. Instead, it seeks to lessen the worst-case scenario, making it appropriate for circumstances where one must protect against every conceivable realization of the uncertain input data.
\par The particular application and its requirements heavily influence whether SO or RO should be used.
The idea of RO, which was initially developed to address USOPs, has created a variety of theories regarding what defines a robust solution to uncertain problems. One well-known idea is minimax robustness, which was first proposed by Soyster \citet{soyster1973convex} and has since been thoroughly investigated in the literature, as can be seen in Ben-Tal et al. \citet{ben2009robust}.
\par The idea of minimax robustness seeks to identify a solution that works in all possible scenarios and provides protection from the worst-case scenario. Finding minimal robust optimal solutions is a difficult two-stage task. There are numerous additional ways to interpret robustness for USOPs. For further details see in (Ben-Tal et al. \citet{ben2009robust}, Ben-Tal and Nemirovski \citet{ben1999robust}, Ben-Tal and Nemirovski \citet{ben1998robust}, Kouvelis and Yu \citet{kouvelis1996robust}, Fischetti and Monaci \citet{fischetti2009light}, Sch$\ddot{\text{o}}$bel \citet{schobel2014generalized}, Erera et al. \citet{erera2009robust}, Liebchen et al. \citet{liebchen2009concept} , Goerigk and Sch$\ddot{\text{o}}$bel \citet{goerigk2015algorithm}).
\par Real-world optimization problems frequently involve multiple objectives and are quite complex. This indicates that these problems frequently require a multiobjective approach in addition to dealing with uncertainty. 
\par Dealing with uncertainties in MOPs, RO is a relatively new area of research. Early research in this field did not focus on the typical concepts of robustness. Instead, most of the work centered around an idea introduced by Branke \citet{branke1998creating}, which was initially meant for USOPs. In this concept, the objective function is replaced by its average value within a specific neighborhood around the chosen point.
\par Based on this concept, Deb and Gupta \citet{deb2013multi} introduced, two approaches to handling uncertainties in MOPs. The objective vector in the first method is changed to a vector that contains the average values of each original component. Therefore, robust solutions to the original problem are efficient solutions to the modified problem. Mean functions are included in the constraints in the second method. This makes sure the objective components do not deviate from their average values by a distance greater than a set limit. Because users can control the desired level of robustness through this threshold, the latter approach is thought to be more practical. There are numerous additional ways to interpret robustness for MOPs. For further details see in (Barricomand Antunes \citet{barrico2006robustness}, Gunawan and Azarm \citet{gunawan2005multi}, Dellnitz and Witting \citet{dellnitz2009computation}, Witting et al.~\citet{witting2013variational}, Kuroiwa and Lee \citet{kuroiwa2012robust}, Fliege and Werner \citet{fliege2014robust}, Yu and Liu \citet{yu2013robust}, Chen et al.~\citet{chen2012including}).
Authors solved the UMOP with the help of the scalarization method (e.g., weighted sum approach and $\epsilon-$constraint approach) via using a general extension of the minimax concept of the robustness from USOP to UMOP. Using RO approach, Ehrgott et al. \citet{ehrgott2014minmax} changed the UMOP into a deterministic MOP with the help of a minimax type robust counterpart and objective-wise worst-case cost type robust counterpart. Ehrgott et al.~\citet{ehrgott2014minmax} solved robust counterparts of UMOPs with the help of scalarization techniques like the weighted sum approach and the $\epsilon-$constraints approach. There are some drawbacks to using scalarization techniques to solve the robust counterparts of UMOPs, such as the pre-specification of weights, restrictions, or function importance that is unknown beforehand.
\par To overcome all these difficulty we solve robust counterpart of UMOP by conjugate gradient descent method. This approach does not make use of predetermined weighting factors or any other kind of ranking or ordering data for the various objective functions. 
Initially the descent methods (or iterative optimization methods) developed for scalar optimization problem. 
\par In order to address the deterministic smooth MOP, a number of descent methods that were initially developed to address SOP have been extended and studied. The concept of iterative methods for solving MOPs was first introduced in  \citet{fliege2000steepest}. Since then, several authors have expanded upon this area, including the development of Newton's method \citet{fliege2009newton}, quasi Newton method \citet{ansary2015modified,povalej2014quasi,lai2020q,mahdavi2020superlinearly,morovati2018quasi,qu2011quasi}, conjugate gradient method \citet{gonccalves2020extension,lucambio2018nonlinear}, projected gradient method \citet{cruz2011convergence,fazzio2019convergence,fukuda2013inexact,fukuda2011convergence,zhao2022linear,drummond2004projected}, and proximal gradient method \citet{bonnel2005proximal,ceng2010hybrid}. Convergence properties are a common characteristic of these methods. Notably, these approaches do not require transforming the problem into a parameterized form, distinguishing them from scalarization techniques \citet{eichfelderadaptive} and heuristic approaches \citet{laumanns2002combining}. It is noticed that these methods are given for deterministic MOP and do not work for uncertain multiobjective optimization.
\par To the best of our knowledge, there is no conjugate gradient method developed for the robust counterpart of UMOP. Although, steepest descent method \citet{kumar2024steepest}, Newton's method \citet{shubham2023newton}, quasi-Newton method \citet{mahato2023quasi}, and modified quasi-Newton method \citet{kumar2024modified} are developed for the robust counterpart of UMOP. 
The steepest descent method \citet{kumar2024steepest}, which uses only the first derivatives (gradients) to determine a search direction, is not always the most efficient. Utilizing higher derivatives can yield algorithms that outperform the steepest descent method. Newton's method \citet{shubham2023newton}, for instance, employs both first and second derivatives and typically performs better than steepest descent, especially when the starting point is close to the local efficient solution. It is also presented that Newton's method has a quadratic rate of convergence. However, for a general nonlinear objective function, convergence to a solution cannot be guaranteed from an arbitrary initial point. If the initial point is not sufficiently close to the solution, the algorithm may not maintain the descent property. A computational drawback of Newton's method \citet{shubham2023newton} is the requirement to evaluate the inverse of the Hessian at each iteration. To avoid computing the Hessian inverse, a quasi-Newton method \citet{mahato2023quasi} is developed to address the solution of the robust counterpart of UMOP. In quasi-Newton method \citet{mahato2023quasi}, the Hessian of each objective function is approximated by a positive definite matrix, which incurs a significant computational cost. To overcome this burden, a modified quasi-Newton method \citet{kumar2024modified} has been developed to solve the robust counterpart of UMOP. In this modified quasi-Newton method, we reduce computational load by generating a common positive definite matrix rather than computing the Hessian for each objective function separately. Briefly, Newton's method \citet{shubham2023newton} is computationally expensive, as it requires both computing and inverting the Hessian matrix, which can be challenging for large-scale problems. Similarly, the quasi-Newton \citet{mahato2023quasi} and modified quasi-Newton \citet{kumar2024modified} methods might fail or be computationally prohibitive for large-scale problems. The conjugate gradient method, however, does not require the computation or inversion of the Hessian matrix, making it more computationally efficient than Newton?s method. Additionally, it eliminates the need to approximate the Hessian inverse for each objective function (or common Hessian approximation). Therefore, we develop a conjugate gradient method for the robust counterpart of UMOP.
\par The following describes how the paper is organized: Important results, fundamental definitions, and theorems essential to our problem are presented in section $\ref{s2}$. The basic concepts of conjugate gradient method are presented in section \ref{secc3}. Using the concepts of nonlinear conjugate gradient method for robust counterpart of UMOP, a nonlinear conjugate algorithm is developed in section \ref{s3}. By using this algorithm, we generate a sequence, and its convergence to a critical point is proven in subsection $\ref{ss3.3}$. In section $\ref{sec5},$ the proposed method is numerically verified using specific test problems. The comparison of the proposed method is also demonstrated with the weighted sum method using performance profiles. In section $\ref{sec6},$ concludes the paper with some remarks on the nonlinear conjugate gradient method.
\section{General concepts in uncertain multiobjective optimization and determinitic multiobjective optimization}\label{s2}
Let us begin with some notations. The set of real numbers denotes as $\mathbb{R}$, the set of non-negative real numbers denotes as $\mathbb{R}_{\geq}=\{x \in \mathbb{R} : x \geq 0\}$, the set of positive real numbers denotes as $\mathbb{R}_{>} = \{x \in \mathbb{R} : x > 0\}$,
$\mathbb{R}^{n} = \mathbb{R} \times \cdots \times \mathbb{R}$ ($n$ times),
$\mathbb{R}^{n}_{\geqq} = \mathbb{R}_{\geq} \times \cdots \times \mathbb{R}_{\geq}$ ($n$ times), $\mathbb{R}^n_{\geq}= \mathbb{R}^{n}_{\geqq}\setminus \{0\},$ $\mathbb{R}^{n}_{>} = \mathbb{R}_{>} \times \cdots \times \mathbb{R}_{>}$ ($n$ times) $\big(\text{or}~ \mathbb{R}^{n}_{>}=int(\mathbb{R}^{n}_{\geqq})\big)$.
In the context of $s, q \in \mathbb{R}^n$, the notation $s \leq/\leqq/< q$ denotes that the component of $s$ is related to the corresponding component of $q$ such that $s_i \leq/\leqq/< q_i$ for all $i = 1, 2, ..., n$.
Similarly, for any $s, q \in
\mathbb{R}^n$: $ s\geqq q\Longleftrightarrow	s - q \in \mathbb{R}^n_{\geqq}$ which is equivalent to $s_i - q_i \geqq 0$ for each $i$; for any $s, q \in
\mathbb{R}^n$: $ s\geq q\Longleftrightarrow	s - q \in \mathbb{R}^n_{\geq}$, which is equivalent to $s_i - q_i \geq 0$ for each $i$; and for any $s, q \in
\mathbb{R}^n:$ $s > q \Longleftrightarrow	 s - q \in \mathbb{R}^n_{>}$, which is equivalent to $s_i - q_i > 0$ for each $i$. It can be noticed that $\mathbb{R}^n_{\geq}=\mathbb{R}^n_{\geqq}$ for $n=1$ i.e., $\mathbb{R}_{\geq}=\mathbb{R}_{\geqq}.$ Throughout the paper, by $A\symbol{92}\{a\},$ we mean set of all elements of $A$ other than $a.$ 
Lastly, we denote the indexed sets as $	\langle p \rangle = \{1, 2, ..., p\}$ and $\langle m \rangle = \{1, 2, ..., m\}$, which contain $p$ and $m$ elements, respectively. Note that, we will use a continuously differentiable function as a CD function. 
\par A MOP can be considered as 
$P:~\displaystyle\min_{x\in D \subset \mathbb{R}^n}F(x),$ where $F:D\to\mathbb{R}^m$ and $F(x)=(F_{1}(x),F_{2}(x),$\ldots$,F_m(x)).$  Now we define Pareto optimal solution (efficient solution) and weak Pareto optimal solution (weak efficient solution) for problem $P.$ A feasible point $x^*\in D$ is considered to be an efficient (weakly efficient) solution for $P$ if there is no another $x\in D$ such that $F(x)\leq F(x^*)$ $\&$  $F(x)\not =F(x^*)$ $\big(F(x)<F(x^*)\big).$ If $x^*$ is an efficient (weakly efficient) solution, then $F(x^*)$ is called non dominated  (weakly non dominated) point, and the set of efficient solution and non dominated point are called efficient set and non dominated set, respectively.  
\par In uncertain multiobjective optimization, we handle input data that are uncertain and affect how we formulate the optimization problem. This uncertainty is encapsulated as a set $U\subset \mathbb{R}^k,$ encompassing all potential scenarios or realizations of the uncertain data. For each realization $\omega$ in the set $U,$ we have a unique MOP 
\begin{equation*}
	P(\omega):~\min_{x\in D}h(x,\omega),
\end{equation*}
where $h:D\times U\to \mathbb{R}^m.$
\begin{rem}
	In this situation, it is crucial to keep in mind that the uncertainty specifically related to the objective function and not the constraints. This is motivated by the fact that a solution to the robust problem is only considered feasible under the usual definition of minimax robustness if it is feasible for every scenario. This type of uncertainty is known as parameter uncertainty.
\end{rem}
\begin{rem}
	When accounting for uncertainty in the variable $D$, a robust solution $x$ to the uncertain problem must adhere to the constraints $x\in \displaystyle\cap_{\omega\in U} D(\omega)$, where $D(\omega)$ represents the feasible set associated with the scenario $\omega.$ This type of uncertainty is known as decision uncertainty.
\end{rem}Throughout this paper, we will take a UMOP with parameter uncertainty.
Considering an uncertainty set $U\subset \mathbb{R}^k$, a feasible set $D\subset \mathbb{R}^n$, and an objective vector-valued function $h:D\times U\to \mathbb{R}^m$, we define an UMOP denoted as $P(U)$, representing the family of $P(\omega)$  such that 
\begin{equation}
	P(U)=\{P(\omega):\omega\in U\}.
\end{equation}
One can observe that for $m=1$, $P(U)$ is an USOP. Additionally, for any given $x\in D$ to $P(U)$, the set of images of $x$ under all scenario is given by $h_U(x)=\{h(x,\omega):\omega\in U\}.$ 
It can be observe that $P(U)$ is a set valued optimization problem. So the concept of optimality of $P(U)$ can be adopted from the set valued optimization. We need a way to compare the sets that represent the various outcomes to find the best solution to a set-based problem. The so-called set approach, as described in Eichfelder and Jahn \citet{eichfelder2011vector}, Ha and Jahn \citet{jahn2011new}, and Kuroiwa \citet{kuroiwa1998natural}, compares sets using order relations. In this context, a set order relation $\preceq$ is employed to compare sets using a specified closed, convex, and pointed solid cone $C\subset\mathbb{R}^m$. The relation $\preceq_C,$ $\preceq_{C\symbol{92}\{0\}},$ and $\preceq_{int\{C\}}$ are defined as follows for every sets $P,~ Q\subset \mathbb{R}^m,$
\begin{eqnarray*}
	P\preceq_C Q \iff P\subset \{Q\}-C,\\
	P\preceq_{C\symbol{92}\{0\}} Q \iff P\subset \{Q\}-C\symbol{92}\{0\},\\\mbox{and}~~~~
	P\preceq_{int\{C\}} Q \iff P\subset \{Q\}-int\{C\}.
\end{eqnarray*}
Equivalently, these set relations can be represented as 
\begin{eqnarray}
	P[\preceq_{int\{C\}}/{\preceq_{C\symbol{92}\{0\}}}/\preceq_C ]\iff \text{for all}~ p\in P,~ \exists ~q\in Q ~\text{such that}~[p</\leq/\leqq q].
\end{eqnarray} 
Throughout this paper, we will take a closed convex pointed solid cone $C=\mathbb{R}^m_\geqq.$ On the behalf of this cone the concept of robust efficiency can be define in following way. \\
\begin{dfn}(Ehrgott et al. \citet{ehrgott2014minmax}) Given an UMOP $P(U),$ a feasible point $x^*\in D$ is said to be 
	\begin{itemize}
		\item robust efficient (RE) if there is no  $x\in D\symbol{92}\{x^*\}$ such that $h(x;U)\subset h(x^*;U)-\mathbb{R}^m_\ge,$
		\item robust weakly efficient (RWE) if there is no  $x\in D\symbol{92}\{x^*\}$ such that $h(x;U)\subset h(x^*;U)-\mathbb{R}^m_>$,
		\item robust strictly efficient (RSE) if there is no  $x\in D\symbol{92}\{x^*\}$ such that $h(x;U)\subset h(x^*;U)-\mathbb{R}^m_\geqq$.
	\end{itemize}
\end{dfn}
The notion of minimax robustness for USOPs was initially proposed by Soyster \citet{soyster1973convex} and studied by Ben-Tal and Nemirovski \citet{ben1998robust}.  In \citet{ehrgott2014minmax}, Ehrgott et al. introduced an extension of the concept of minimax robustness for USOP to UMOP.
\par Given an UMOP $P(U),$ the concept of minimax robustness is presented by Ehrgott et al. \citet{ehrgott2014minmax}. The motive of this robustness concept is to search for solutions which minimize the worst case that occur, i.e., solution to the problem
\begin{eqnarray}\label{minimax}
	\displaystyle 	\min_{x\in D}\sup_{\omega\in U}h(x,\omega),~\text{where}~ h:D\times U\to \mathbb{R}^m.
\end{eqnarray}	
The authors emphasized the lack of clarity in defining a worst case, attributed to the absence of a total order on $\mathbb{R}^m$. Consequently, a direct extension of the concept of minimax robustness to UMOP was not possible. Therefore, Ehrgott et al. \citet{ehrgott2014minmax} presented an extension of minimax robustness to MOPs, namely the concept of robust efficiency. This concept is presented above in the form of RE solution, RWE solution, and RSE solution. 
\par Ehrgott et al. \citet{ehrgott2014minmax} developed algorithms for computing minimax RE solution. Firstly, a well known weighted sum scalarization approach for deterministic MOP is modified to discover minimax RE solutions to UMOPs. 
\par To find the solution of $P(U)$, a deterministic SOP is defined for a given $w\in\mathbb{R}^m_{\geq}$ such that
\begin{eqnarray*}
	(WP)_{P(U)}(w):~~~\displaystyle \min_{x\in D}\sup_{\omega\in U}\sum^{m}_{i=1}w_ih_i(x,\omega).
\end{eqnarray*}
It is also given that the solution of $(WP)_{P(U)}(w)$ will be the solution for $P(U).$	
\begin{thm}( Ehrgott et al. \citet{ehrgott2014minmax})
	Given an UMOP $P(U),$ the following conditions hold.
	\begin{enumerate}[(i)]
		\item If $x'$ is the unique optimal solution to $(WP)_{P(U)}(w)$ for some $w\in\mathbb{R}^m_{\geq}$, then $x'$ is minimax RSE for $P(U).$
		\item If $x'$ is an optimal solution to $(WP)_{P(U)}(w)$ for some $w\in\mathbb{R}^m_{>}$ and $\displaystyle\max_{\omega\in U}\sum^{m}_{i=1}w_ih_i(x,\omega)$ exist for all $x\in D$, then $x'$ is minimax RE for $P(U).$
		\item If $x'$ is an optimal solution to $(WP)_{P(U)}(w)$ for some $w\in\mathbb{R}^m_{\geq}$ and $\displaystyle\max_{\omega\in U}\sum^{m}_{i=1}w_ih_i(x,\omega)$ exist for all $x\in D$, then $x'$ is minimax RE for $P(U).$
	\end{enumerate}
\end{thm}
Ehrgott et al. \citet{ehrgott2014minmax} presented second approach for calculating the minimax RE solution for $P(U),$ which is an extension of  $\epsilon-$constraints scalarization method for deterministic multiobjective optimization problem. With this approach, to find the minimax RE solution for $P(U),$ we solve following deterministic single objective optimization problem for $\epsilon=(\epsilon_1,\epsilon_2,\ldots,\epsilon_m)\in \mathbb{R}^m_\geq.$  
\begin{eqnarray}
	\epsilon_{P(U)}:~\displaystyle \min_{x\in D}\sup_{\omega\in U}h_l(x,\omega) \nonumber\\
	~~~~~~~~~~~~~~~~~~~~~~~s.t.~~~h_j(x,\omega)\leq\epsilon_j,~\text{for all}~j\not=l,~\omega\in U \nonumber.
\end{eqnarray}	
In the next results it can be observe that the optimal solutions to $\epsilon_{P(U)}$ will be the solution of $P(U).$
\begin{thm}\label{theoremowc}( Ehrgott et al. \citet{ehrgott2014minmax}). Given an UMOP $P(U)$, the following conditions hold.
	\begin{enumerate}[(i)]
		\item If $x'$ is the unique optimal solution to $\epsilon_{P(U)}$ for some $\epsilon\in \mathbb{R}^m$ and for some $l\in\{1,2,\ldots,m\},$
		then $x'$ is minimax RSE for $P(U).$
		\item If $x'$ is the optimal solution to $\epsilon_{P(U)}$ for some $\epsilon\in \mathbb{R}^m$ and for some $l\in\{1,2,\ldots,m\},$
		$\displaystyle\max_{\omega\in U}h_i(x,\omega)$ exists for all $x\in D$, then $x$ is minimax RWE for $P(U).$
	\end{enumerate}
\end{thm}
One more approach is presented by Ehrgott et al. \citet{ehrgott2014minmax} for finding the minimax robust efficient solution to $P(U).$ In this approach a objective wise worst case method is considered. This approach consists the same deterministic multiobjective optimization problem as Kuroiwa and Lee \citet{kuroiwa2012robust} considered for their concepts of mutiobjective robustness, namely    
\begin{eqnarray}\label{owc}
	OWC_{P(U)}:~~\min_{x\in D} h^{owc}_U(x),~\text{where}~ h^{owc}_U(x)=(\sup_{\omega\in U}h_1(x,\omega),\sup_{\omega\in U}h_2(x,\omega),\ldots,\sup_{\omega\in U}h_m(x,\omega))^T.
\end{eqnarray}
Ehrgott et al. \citet{ehrgott2014minmax} assert that solving the $OWC_{P(U)}$ problem will yield the solution to the $P(U)$ problem.
\begin{thm}\label{t1}( Ehrgott et al. \citet{ehrgott2014minmax}) Let $P(U)$ be an UMOP. Then,
	\begin{enumerate}[(a)]
		\item if $x^*\in X$ is a strictly efficient solution to $OWC_{P(U)},$ then $x^*$ is RSE solution for $P(U)$.
		\item If $\displaystyle\max_{\omega_i\in U} h_j(x,\omega)$  exist for all $j\in	\langle m \rangle$ and $x \in D,$  $x^*$
		is a weakly efficient solution to $OWC_{P(U)},$ then $x^*$ is RWE solution for $P(U)$.
	\end{enumerate}
\end{thm}
In \citet{ehrgott2014minmax}, weighted sum scalarization method and  $\epsilon-$constraints scalarization method are also developed for $OWC_{P(U)}$ to find the solution of $P(U).$ It is also proved that the solution obtained by both scalarization method for $OWC_{P(U)}$ will be the solution of $P(U).$
\begin{rem} Obviously, computing $h^{owc}_U(x)$ for any given $x$ is much easier than solving a MOP $\max\{h(x,\omega):\omega\in U\}.$ Since $OWC_{P(U)}$ contains only $m$ deterministic SOPs, then $OWC_{P(U)}$ is a deterministic MOP.
\end{rem}
\begin{rem}
	Theorem \ref{theoremowc}, shows that the solution of $OWC_{P(U)}$ will be the solution for $P(U).$ Also, the problem $OWC_{P(U)}$ is solved by the scalarization method (e.g., weighted sum method and $\epsilon-$constraint method). By pre-selecting some parameters and reformulating them as deterministic SOP, scalarization methods, which are based on the scalarization technique, compute the efficient or weak efficient solution. This method has a drawback in that the parameter selection may result in unbounded (i.e., no solution exists) scalarized problem even when the original deterministic MOP, $OWC_{P(U)}$ has solutions. Another drawback of this strategy is that the parameters are not predetermined, so it is up to the modeler and the decision-maker to make those decisions. To remove all these difficulties, we will solve $OWC_{P(U)}$ with the help of the NCGM, and there is no parameter information is needed for this method.
\end{rem}
\begin{rem}
	Conjugate gradient method is an iterative method (numerical optimization method). This approach does not make use of predetermined weighting factors or any other kind of ranking or ordering data for the various objective functions.
\end{rem}
We will develop a conjugate gradient method for $OWC_{P(U)}$ to find the solution of $P(U).$ We consider that $U$ is a finite uncertainty set and $D=\mathbb{R}^n$. In case of finite uncertainty, if we assume $U^F=\big\{\omega_i:i\in \langle p \rangle=\{1,2,\ldots,p\}\big\}$ and $h_j(x,\omega_i)$ is CD for each $x$ and $\omega_i$ then the problem $P(U)$ and $OWC_{P(U)}$ can be considered as $P(U^F)$ and $OWC_{P(U^F)},$ respectively, which are defined as follows:
\begin{equation}\label{main uncertain mop}
	P(U^F)=\{P(\omega_i):\omega_i\in U^F\}
\end{equation}
and
\begin{equation*}\label{main owc}
	OWC_{P(U^F)}:~~\displaystyle\min_{x\in \mathbb{R}^n} h^{owc}_{U^{F}}(x),
\end{equation*}
where $h^{owc}_{U^{F}}(x)=(\displaystyle\max_{\omega_i\in U^F}h_1(x,\omega_i),\displaystyle\max_{\omega_i\in U^F}h_2(x,\omega_i),\ldots,\displaystyle\max_{\omega_i\in U^F}h_m(x,\omega_i))^T.$\\
Also, we can write $OWC_{P(U^F)}$ as follows:
\begin{eqnarray}\label{main owc}
	OWC_{P(U^F)}:~~\min_{x\in \mathbb{R}^n} \digamma(x),
\end{eqnarray}
where $\digamma(x)=(\digamma_1(x),\digamma_2(x),\ldots,\digamma_m(x))^T~\text{and}~\digamma_j(x)=\displaystyle\max_{\omega_i\in U^F}h_j(x,\omega_i),~	j\in\langle m \rangle.$\\
Since $U^F$ is a finite set, and $h_j(x, \omega_i)$ is a CD function for each $x$ and $\omega_i$, $\displaystyle\sup_{\omega_i\in U^F}h_j(x,\omega_i)=\displaystyle\max_{\omega_i\in U^F}h_j(x,\omega_i)$ due to the compactness of $U^F.$
\par \par In this paper, we will develop a NCGM for $OWC_{P(U^F)}$ to find the solution of $P(U^F).$ We are mentioning NCGM here because we are considering nonlinear objective functions in $P(U^F).$ One can observe that $OWC_{P(U^F)}$ is a non-smooth deterministic MOP.
Now, we define necessary condition of Pareto optimality for $OWC_{P(U^F)}.$
\begin{dfn}\label{criticalpoint} A point $x^*\in D$ is said to be critical point for $OWC_{P(U^F)}$ if $$I(\partial\digamma(x^*))\cap(-\mathbb{R}^{n}_{>})= \emptyset,$$  where
	$\partial\digamma(x^*)=Conv\{\displaystyle \cup_{j\in 	\langle m\rangle} \partial \digamma_j(x^*)\}.$
	It is noticed that $x^*$ is a critical point for $OWC_{P(U^F)}$ if and only there is no $v\in \mathbb{R}^n$ such that $\nabla h_j(x^*,\omega_i)^Tv<0$, for all $i\in I_j(x),$ $j\in 	\langle m\rangle$, where $I_j(x)=\{i\in	\langle p \rangle:h_{j}(x^*,\omega_i)=\digamma_j(x^*)\}.$ \end{dfn}
\begin{note} By definition of critical point for $OWC_{P(U^F)}$ and Theorem \ref{t1}, we can conclude that critical point for $OWC_{P(U^F)}$ will be the robust critical point for $P(U^F).$ Throught the paper, we will use critical point for $OWC_{P(U^F)}$ as the robust critical point for $P(U^F).$
\end{note}
\begin{dfn}\citet{kumar2024steepest}
	In $OWC_{P(U^F)},$ a vector $v$ is said to be descent direction for $\digamma$  at $x$ if
	$\nabla h_j(x,\omega_i)^Tv< 0, ~\forall~j \in 	\langle m\rangle$ and $i\in I_j(x).$
	Also, $v$ is descent direction for $\digamma(x)$ $\Longleftrightarrow$ there exists $\epsilon >0$ such that $\digamma_j(x+\alpha v)< \digamma_j(x), ~\forall~ j \in 	\langle m\rangle~ \text{and} ~\alpha \in(0,\epsilon].$
\end{dfn}
\begin{thm}\label{directional derivative}\citet{dhara2011optimality} Let $\digamma_j: \mathbb{R}^n \to \mathbb{R}^n$ denote a function such that $\digamma_j(x)=\displaystyle\max_{i\in	\langle p \rangle}h_j(x,\omega_i).$ Then:
	\begin{enumerate}[(i)]
		\item In the direction $v$, the directional derivative of $\digamma_j$ at $x$ is given by $H'_j(x,v)=\displaystyle\max_{i\in I_j(x)}\nabla h_j(x,\omega_i)^Tv,$ where $I_j(x) = \{i \in 	\langle p \rangle : h_j(x, \omega_i) = \digamma_j(x)\}$.
		\item\label{ddii} The subdifferential of $\digamma_j$ is  $\partial \digamma_j(x)= Conv \bigg(\bigcup_{i\in I_j(x)}\partial h_j(x,\omega_i)\bigg).$ Also, $x^*=\displaystyle \underset{x\in \mathbb{R}^n}{\mathrm{argmin}}\digamma_j(x)$ $\Longleftrightarrow$ $0\in\partial \digamma_j(x^*)).$ 
	\end{enumerate}
\end{thm}
To demonstrate the necessary and sufficient condition for the Pareto optimality of $OWC_{P(U^F)}$, we first provide the following Lemma \ref{lm1}.
\begin{lem}\label{lm1}\citet{kumar2024steepest}
	The point $x^*$ is a critical point for $\digamma$ $\Longleftrightarrow$ $0\in Conv\{\displaystyle \cup_{j\in 	\langle m\rangle} \partial \digamma_j(x^*)\}$.
\end{lem}
\begin{thm}\label{nce}
	If  $h_j(x,\omega_i)$ is CD and convex for each $j\in 	\langle m\rangle$ and $\omega_i\in U$, then  $x^*\in \mathbb{R}^n$ is a weakly Pareto optimal solution for $OWC_{P(U^F)}$ if and only if
	$$	0 \in conv \left( \displaystyle \cup_{j=1}^{m}\partial \digamma_j(x^*) \right).$$
\end{thm}
\begin{proof}
	See Theorem 2.3 in \citet{kumar2024steepest}.
\end{proof}
\section{Nonlinear conjugate gradient method for multiobjective optimization}\label{secc3}
As we mentioned in the introduction section, we develop conjugate gradient method for uncertain multiobjective optimization problem. So first we introduce conjugate gradient method for deterministic optimization problem. Conjugate gradient method originally proposed by R. Fletcher and C.M. Reeves in 1964 for unconstrained scalar  optimization problem. L. R. Lucambio Perez and L. F. Prudente extended the idea of conjugate gradient method for scalar optimization to vector optimization, for further details see \citet{lucambio2018nonlinear}. In \citet{gonccalves2020extension,lucambio2018nonlinear}, authors developed conjugate gradient algorithm for smooth multiobjective (vector) optimization problem $\displaystyle\min_{x\in \mathbb{R}^n}f(x),$ where $f(x)=(f_{1}(x),f_{2}(x),\ldots,f_m(x)).$ In this algorithm, authors generated a sequence of iterates according to $x^{k+1}=x^k+\alpha_kv^k,~k\geq0,$ where $v^k\in\mathbb{R}^n$ is the line search direction and $\alpha_k$ is the step size. The direction $v^k$ is defined as $v^k=s(x^k)+\beta_kv^{k-1},$ $k\geq 1,$ where $\beta_k$ is a scalar algorithmic parameter, $v^0=s(x^0),$ and $s(x^k)={\arg\min}_{d\in\mathbb{R}^n}\{f(x^k,d)+\frac{\|d\|^2}{2}\}.$ For non quadratic functions different choice of $\beta_k$ is considered. Some remarkable choice of $\beta_k$ are as follows:\\
\textbf{Fletcher-Reeves (FR):} $\beta_k=\frac{f(x^k,s(x^k))}{f(x^{k-1},v(x^{k-1}))},$ where $f(x^k,s(x^k))=\displaystyle\max_{j=1,\ldots,m}\{\nabla f_j(x)^Ts(x^k)\}.$
\textbf{Conjugate descent (CD):} $\beta_k=\frac{f(x^k,s(x^k))}{f(x^{k-1},v^{k-1})},$ where $f(x^k,s(x^k))=\displaystyle\max_{j=1,\ldots,m}\{\nabla f_j(x)^Ts(x^k)\}.$
\textbf{Dai-Yuan (DY):} $\beta_k=\frac{-f(x^k,s(x^k))}{f(x^{k},v^{k-1})-f(x^{k-1},v^{k-1})},$ where $f(x^k,s(x^k))=\displaystyle\max_{j=1,\ldots,m}\{\nabla f_j(x)^Ts(x^k)\}.$\vspace{0.2cm}
\textbf{Polak-Ribiere-Polak (PRP):}  $\beta_k=\frac{-f(x^k,s(x^k))+f(x^{k-1},s(x^k))}{-f(x^{k-1},v(x^{k-1}))}.$\vspace{0.2cm}\\
\textbf{Hestenes--Stiefel (HS):}  $\beta_k=\frac{-f(x^k,s(x^k))+f(x^{k-1},s(x^k))}{f(x^k,v^{k-1}))-f(x^{k-1},v^{k-1})}.$\\
\par Next, we formulate certain theoretical concepts associated with the nonlinear conjugate gradient method for the $OWC_{P(U^F)}$ problem.
\section{Nonlinear conjugate gradient method for $OWC_{P(U^F)}$}\label{s3}
Let  $ \digamma:\mathbb{R}^n \to \mathbb{R}^m$ be such that $\digamma(x)=(\digamma_1(x),\digamma_2(x),\ldots,\digamma_m(x))$, where $\digamma_j(x)= \displaystyle \max_{i\in 	\langle p \rangle} h_j(x,\omega_i),\\~ j\in \langle m \rangle$ and $h_j:\mathbb{R}^n \times U^F \to \mathbb{R}$ is a CD and convex function for each $x\in \mathbb{R}^n$ and $\omega_i\in U^F$. Then, $v\in\mathbb{R}^n$ is said to be descent direction of $\digamma$ at $x$ if and only if $ \nabla h_j(x,\omega_i)^Tv <0,~ \forall ~i \in I_j(x)$ and $j \in 	\langle m\rangle $. Now we define a function $h:\mathbb{R}^n \times\mathbb{R}^n\to \mathbb{R}$ such that 
\begin{equation}\label{4.1}
	h(x,v)=\displaystyle  \max_{j\in 	\langle m\rangle}	\displaystyle  \max_{i\in 	\langle p \rangle}\{h_j(x,\omega_i)+\nabla h_j(x,\omega_i)^Tv-\digamma_j(x)\}
\end{equation}
\begin{lem}\label{lem1}
	Let $\digamma:\mathbb{R}^n\to\mathbb{R}^m$ be a function such that $\digamma(x)=(\digamma_1(x),\digamma_2(x),\ldots,\digamma_m(x))$, where $\digamma_j(x)= \displaystyle \max_{i\in 	\langle p \rangle} h_j(x,\omega_i),~ j\in 	\langle m\rangle$ and $h_j:\mathbb{R}^n \times U^F \to \mathbb{R}$ is a CD for each $i$ and $j$. Then the following statements are equivalent:
	\begin{enumerate}[(1)]
		\item\label{l1} $h(x,v)<0.$
		\item\label{l2} $\nabla h_j(x,\omega_i)^Tv<0,$ for all $i$, $j$ provided $h_j(x,\omega_i)-\digamma_j(x)=0,$ for all $i\not\in I_j(x).$
		\item\label{l3} $v$ is a descent direction. 
	\end{enumerate}
\end{lem}
\begin{proof}
	$(1)\implies (2):$\\
	By item (\ref{l1}), $h(x,v)<0.$ Then, 
	\begin{align*}
		\displaystyle  \max_{j\in 	\langle m\rangle}	\displaystyle  \max_{i\in 	\langle p \rangle}\{h_j(x,\omega_i)+\nabla h_j(x,\omega_i)^Tv-\digamma_j(x)\}<0,\\
		\implies h_j(x,\omega_i)+\nabla h_j(x,\omega_i)^Tv-\digamma_j(x)<0, ~\text{for all}~ i\in\langle p\rangle,~j\in\langle m\rangle.
	\end{align*}
	Since we know for each $i\in I_j(x)$, $h_j(x,\omega_i)=\digamma_j(x)$ and also if we consider $h_j(x,\omega_i)-\digamma_j(x)=0$ for each $i\notin I_j(x)$ then item (\ref{l2}) follows.\\ 
	$(2)\implies (3):$\\
	By item (\ref{l2}), $\nabla h_j(x,\omega_i)^Tv<0,$ for all $i$, $j$ then $v$ is a descent direction for $\digamma$ at $x$ hence the item (\ref{l3}).\\ 
	$(3)\implies (1):$\\
	By item (\ref{l3}), $v$ is a descent direction for $\digamma$ at $x$ provided $h_j(x,\omega_i)-\digamma_j(x)=0,$ for all $i\not\in I_j(x),$ i.e., $\nabla h_j(x,\omega_i)^Tv<0,$ for all $i$, $j$. Then, 
	\begin{align*}
		h(x,v)&=\displaystyle  \max_{j\in 	\langle m\rangle}	\displaystyle  \max_{i\in 	\langle p \rangle}\{h_j(x,\omega_i)+\nabla h_j(x,\omega_i)^Tv-\digamma_j(x)\}\\
		&\leq\displaystyle  \max_{j\in 	\langle m\rangle}\{	\displaystyle  \max_{i\in 	\langle p \rangle}h_j(x,\omega_i)+	\displaystyle  \max_{i\in 	\langle p \rangle}\nabla h_j(x,\omega_i)^Tv-\digamma_j(x)\}\\
		&\leq\displaystyle  \max_{j\in 	\langle m\rangle}	\displaystyle  \max_{i\in 	\langle p \rangle}\nabla h_j(x,\omega_i)^Tv.
	\end{align*}
	Thus, $h(x,v)\leq\displaystyle  \max_{j\in 	\langle m\rangle}	\displaystyle  \max_{i\in 	\langle p \rangle}\nabla h_j(x,\omega_i)^Tv<0.$ Hence item (\ref{l1}).
\end{proof}
%
In the next subsection, we construct a subproblem, and solve it to find the descent direction for $OWC_{P(U^F)}$.
\subsection{A subproblem to find the descent direction for $OWC_{P(U^F)}$}\label{ss3.1}
The notion of steepest descent descent direction for $OWC_{P(U^F)}$ at $x$ is given by 
\begin{eqnarray}
	s(x)&=&\underset{v\in \mathbb{R}^n}{\mathrm{argmin}}  \left( h(x,v)+\tfrac{1}{2}\|v\|^2 \right)\label{ov}\\
	\text{and} ~~~~~T(x)&=&  h(x,v(x))+\frac{1}{2}\|v(x)\|^2,\label{os}
\end{eqnarray}$\text{ where } h(x,v)= \displaystyle \max_{j\in	\langle m\rangle}\max_{i\in 	\langle p \rangle}  \{ h_j(x,\omega_i)+\nabla h_j(x,\omega_i)^Tv-\digamma_j(x)\}.$
To find $s(x)$, we need to solve the subproblem  \begin{equation}\label{eq31}\displaystyle\min_{v\in\mathbb{R}^n} \left( h(x,v)+\tfrac{1}{2}\|v\|^2 \right).
\end{equation}
Note that, as the maximum of the max-linear function $h(x,v)$ is convex, the objective function in the subproblem \eqref{eq31} is strongly convex. Therefore, the solution to this subproblem \eqref{eq31} is a unique solution and serves as the descent direction for $OWC_{P(U^F)}$ (see Theorem \ref {Thm 3.2}).
Equivalently, subproblem \eqref{eq31} can be written as
\begin{eqnarray*}
	P(x):  \min_{v\in\mathbb{R}^n,t\in\mathbb{R}}~t+\frac{1}{2}\|v\|^2&\\
	\text{s.t.}~~ h_j(x,\omega_i)+\nabla h_j(x,\omega_i)^Tv -\digamma_j(x)\leq t,&\forall ~i\in 	\langle p \rangle,~j \in 	\langle m\rangle.
\end{eqnarray*}
One can observe that $P(x)$ satisfies Slater's constraint qualification since the inequalities in $P(x)$ are strict for $t=1$ and $v=(0,0,\ldots,0)\in \mathbb{R}^n$. Hence, there exists $	\langle m\rangle\in\mathbb{R}^{m\times p}_{+}$ such that the following KKT optimality conditions hold:
\begin{eqnarray}
	\sum\limits_{j\in 	\langle m\rangle}  \sum\limits_{i\in 	\langle p \rangle} 	\lambda_{ij} &=&1,\label{kkt1}\\
	v +\sum\limits_{j\in 	\langle m\rangle} \sum\limits_{i\in 	\langle p \rangle} 	\lambda_{ij} \nabla h_j(x,\omega_i)&=&0,\label{kkt2}\\
	h_j(x,\omega_i)+\nabla h_j(x,\omega_i)^Tv -\digamma_j(x) -t&\leq& 0,~ \forall~ i\in 	\langle p \rangle,~j\in	\langle m\rangle, \label{kkt3}\\
	\lambda_{ij}\geq 0,~~~	\lambda_{ij}\left(h_j(x,\omega_i)+\nabla h_j(x,\omega_i)^Tv  -\digamma_j(x)-t\right)&=&0, ~~\forall~ i\in 	\langle p \rangle,~ j\in	\langle m\rangle.\label{kkt4}
\end{eqnarray}
Problem $(P(x))$ has a unique solution, $(s(x),T(x)).$ As this is a convex problem and has Slater point, there exists $\lambda=\lambda(x)=\lambda_{ij}(x),$ together $v=s(x),$ $t=T(x),$ satisfies the conditions (\ref{kkt1}), (\ref{kkt2}), (\ref{kkt3}), and (\ref{kkt4}). Consequently, by (\ref{kkt2}), we get
\begin{equation}\label{4***}
	s(x) = -	\sum\limits_{j\in 	\langle m\rangle}  \sum\limits_{i\in  	\langle p\rangle } 	\lambda_{ij} \nabla h_j(x,\omega_i).
\end{equation} 
Further, $T(x)$ is  the  optimal value of the subproblem \eqref{eq31} which is given as follows
\begin{equation}
	T(x)= h(x,s(x))+\frac{1}{2}\|s(x)\|^2.
\end{equation}
\qed\\
\par Now, let us examine certain properties of the function $T$ and explore its connection with $s(x)$ and the criticality of $x.$ 
\begin{thm}\label{Thm 3.2}\rm Let $s(x)$ and $T(x)$ are defined in Eqn. $(\ref{ov})$ and Eqn. $(\ref{os})$, respectively. Then the following results hold:
	\begin{enumerate}
		\item\label{th0} $ s(x)$ is bounded on compact subset   $C$ of  $\mathbb{R}^n$ and $T(x)\leq0$.
		\item The following conditions are equivalent :
		\begin{enumerate}
			\item The point $x$ is not a critical point.
			\item\label{theta<0} $T(x)<0.$
			\item $s(x)\not =0.$
			\item $s(x)$ is a descent direction for $\digamma$ at $x$.
		\end{enumerate}
	\end{enumerate}
	In particular $x$ is critical point if and only if $T(x)=0.$
\end{thm}
\begin{proof} 
1. Let $C$ be a compact subset of $\mathbb{R}^n$. Since $h_j(x,\omega_i)$ is a convex and differentiable (CD) function for each $j\in \langle m\rangle$ and $i\in \langle p \rangle$, $h_j(x,\omega_i)$ is bounded on every compact set. So, for all $x\in C$, $i\in	\langle p \rangle$ and $j \in 	\langle m\rangle,$ by Eqn. $(\ref{4***}),$ $s(x)$ is bounded on the compact set $C$. Since $t=1$ and $s=\bar0=(0,0,\ldots,0)\in \mathbb{R}^n$ lies in the feasible region then we have\\ 	
$~~~~~~~T(x) \leq \displaystyle \max_{j\in 	\langle m\rangle}\{\max_{i\in	\langle p\rangle}h_j(x,\xi_i) + \max_{i\in	\langle p\rangle}\nabla h_j(x,\xi_i)^T\bar 0 +\max_{i\in \bar\Lambda} \frac{1}{2}\|\bar 0\|^2-\max_{i\in j\in 	\langle p\rangle}\digamma_j(x)\}=0.$\\
Hence $T(x)\leq0$.\\
$2: (a)\implies (b)$	Since $x$ is not a critical point, there exists $\bar v$ such that
\begin{center}
	$\nabla h_j(x,\omega_i)^T \bar v < 0, ~\forall ~j \in 	\langle m\rangle$ and $i \in I_j(x)$.
\end{center}
Since $T(x)$ is the optimal value for the subproblem \eqref{eq31}, for all $\delta >0,$ we have
\begin{align}
	T(x) &\leq \displaystyle \max_{j\in 	\langle m\rangle}\max_{i\in 	\langle p \rangle} \{h_j(x,\omega)+\nabla h_j(x,\omega_i)^T(\delta \bar v)-\digamma_j(x)\}+\frac{1}{2}\|\delta \bar v\|^2\nonumber\\
	&\leq  \displaystyle \max_{j\in 	\langle m\rangle}\{\max_{i\in 	\langle p \rangle} h_j(x,\omega_i)+ \max_{i\in 	\langle p \rangle}\nabla h_j(x,\omega_i)^T (\delta\bar v)-\max_{i\in 	\langle p \rangle} \digamma_j(x)\}+\frac{1}{2}\| \delta\bar v\|^2 \nonumber\\
	&=\displaystyle \max_{j\in 	\langle m\rangle}\max_{i\in 	\langle p \rangle}\nabla h_j(x,\omega_i)^T (\delta\bar v)+\frac{1}{2}\| \delta\bar v\|^2 \nonumber\\
	&=\delta\big(\displaystyle \max_{j\in 	\langle m\rangle}\max_{i\in 	\langle p \rangle}\nabla h_j(x,\omega_i)^T (\bar v)+\frac{1}{2}\delta\| \bar v\|^2 \big).\nonumber
\end{align}
For small enough $\delta>0$, the right hand side of above inequality will be negative because of $\nabla h_j(x,\omega_i)^T\bar v < 0,$ for all $j\in	\langle m\rangle$ and $i\in I_j(x)$.  Thus, $T(x)
< 0$.\\ 
$2: (b)\implies (c)$  Since $T(x)$ is the optimal value of the subproblem (\ref{eq31}), and from $(b)$, it is negative, so  we will get $s(x) \not =0$. If $s(x)  = 0.$ then $T(x)$ will be zero, which is not possible from $(b)$. Hence, if $T(x) < 0,$ then $s(x) \not = 0$.\\ 
$2: (c)\implies (d)$  Let $s(x) \not = 0$. Then, $T (x) \not =0$. Since $T(x) \leq 0$, $T (x)<0$. Thus,
\allowdisplaybreaks
\begin{align*}
	&T(x) = \displaystyle \max_{j\in 	\langle m\rangle}\max_{i\in	\langle p \rangle}\{h_j(x,\omega_i)+\nabla h_j(x,\omega_i)^Ts(x)-\digamma_{j}(x)\} +\frac{1}{2}\|s(x)\|^2 < 0\\
	&\implies \displaystyle \max_{j\in 	\langle m\rangle}\{\max_{i\in	\langle p \rangle} h_j(x,\omega_i)+\max_{i\in 	\langle p \rangle}\nabla h_j(x,\omega_i)^Ts(x)-\max_{i\in \bar	\langle m\rangle}\digamma_j(x)\}+\frac{1}{2}\|s(x)\|^2 < 0\\
	& \implies \max_{j\in 	\langle m\rangle}\max_{i\in 	\langle p \rangle}\nabla h_j(x,\omega_i)^Ts(x)< 0\\
	&\implies\nabla h_j(x,\omega_i)^Ts(x)< 0, ~\forall ~ j\in 	\langle m\rangle ~\text{and}~ i\in 	\langle p \rangle\\
	&\implies\nabla h_j(x,\omega_i)^Ts(x)< 0, ~\forall ~ j\in 	\langle m\rangle ~\text{and}~ i\in I_j(x)\\
	&\implies s(x)~ \text{is a descent direction for}~ \digamma~ \text{at} ~x.
\end{align*}
$2: (d)\implies (a)$
Since $s(x)$ is a descent direction for $\digamma$ at $x$, we have \begin{equation*}
	\nabla h_j(x,\omega_i)^Ts(x)< 0, ~\forall ~ j\in 	\langle m\rangle ~\text{and}~ i\in I_j(x).
\end{equation*}
Then, $x$ is not a critical point. It follows (a).\\
Also,  if $T(x)<0$, then $s(x) \not= 0$. Moreover, if $s(x) \not= 0$, then $T(x)<0$. Thus, $x$ is  critical point if and  only if $T(x) =0$.
\end{proof}
To write the conjugate gradient algorithm for $OWC_{P(U^F)}$, we next need the step length. In the next Subsection \ref{ss3.2}, we give a rule of choosing the step length.
\subsection{Armijo type inexact line search for conjugate gradient method }\label{ss3.2}
In this subsection in Theorem \ref{linesearch}, inexact line search technique is developed to find a suitable step length size that ensures sufficient decrease in each objective function.
\begin{thm}\label{linesearch}
Assume that $x$ is not critical point of $\digamma$. Then for any $\beta>0$  and $\epsilon\in(0,1]$ there exists an $\alpha\in[0,\epsilon]$ such that
\begin{equation*}\label{18''}
	\digamma_j(x+\alpha v)\leq  \digamma_j(x)+  \alpha \beta h(x,v),
\end{equation*}
where $h(x,v)$ is given in Eqn. (\ref{4.1}).
\end{thm}
\begin{proof}
Since $t=0,v=(0,0,0,\ldots,0)\in \mathbb{R}^n$ is feasible for $P(x)$, then
\begin{center}
	$T(x)	= \displaystyle \max_{j\in	\langle m\rangle}\max_{i\in 	\langle p \rangle}  \{h_j(x,\omega_i)+\nabla h_j(x,\omega_i)^Tv-\digamma_j(x)\} +\frac{1}{2} \|v\|^2\leq0,$
\end{center}
which gives
\begin{equation*}
	\displaystyle \max_{j\in	\langle m\rangle}\max_{i\in 	\langle p \rangle}  \nabla h_j(x,\omega_i)^Tv \leq-\frac{1}{2} \|v\|^2.
\end{equation*}
This implies
\begin{equation}\label{7}
	\displaystyle\max_{i\in 	\langle p \rangle}  \nabla h_j(x,\omega_i)^Tv \leq-\frac{1}{2} \|v\|^2,~ ~  \forall ~j\in 	\langle m\rangle.
\end{equation}
To find the step length size rule we define an auxiliary function
\begin{equation}\label{8}
	\digamma^*_j(x,v)= \max_{i\in	\langle p\rangle}\{ h_j(x,\omega_i) + \nabla h_j(x,\omega_i)^Tv \}-\digamma_j(x),~ j\in 	\langle m\rangle.
\end{equation}
Note that Eqn. (\ref{8}) implies
\begin{equation}\label{8**}
	\digamma^*_j(x,v) \leq  \displaystyle\max_{i\in	\langle p \rangle}\nabla h_j(x,\omega_i)^Tv,~~ j\in \langle m\rangle,
\end{equation}
which implies 
\begin{align*}\label{8**}
	\digamma^*_j(x,v) &\leq \displaystyle\max_{j\in {\Lambda}} \displaystyle\max_{i\in 	\langle p \rangle}\{h_j(x,\omega_i)+\nabla h_j(x,\omega_i)^Tv-\digamma_j(x)\}\nonumber\\
	&=h(x,v).
\end{align*}
Therefore, \begin{equation}\label{g}
	\digamma^*_j(x,v) \leq h(x,v).
\end{equation}
Since $x$ is not a critical point, $h(x,v)<0.$
Then, we get
\begin{equation}\label{11***}
	\digamma^*_j(x,v) < 0.
\end{equation}
Since $h_j(x,\omega_i)$ is convex and a CD function, it is thus upper uniformly differentiable (refer to Definition 2.1 on page 159 in \citet{bazaraa1982algorithm}). Then there exists $k^j_i$ such that
\begin{center}
	$	h_j(x+v,\omega_i) \leq h_j(x,\omega_i)+\nabla h_j(x,\omega_i)^Tv + \frac{1}{2} k^j_i \|v\|^2,~~\text{for each}~ i\in 	\langle p \rangle ~\text{and} ~	j\in\langle m \rangle.$
\end{center}
Also,
\begin{align*}
	h_j(x+v,\omega_i) &\leq h_j(x,\omega_i)+\nabla h_j(x,\omega_i)^Tv + \frac{1}{2} k^j_i \|v\|^2 \\ &\leq \max_{i\in 	\langle p \rangle} \{ h_j(x,\omega_i)+\nabla h_j(x,\omega_i)^Tv \}+ \frac{1}{2} K \|v\|^2,
\end{align*}
where $K^j= \displaystyle \max_{i\in 	\langle p \rangle}{k^j_i}$. Now, from the Eqn. (\ref{8}), we have
\begin{equation*}
	h_j(x+v,\omega_i) \leq 	\digamma^*_j(x,v) + \digamma_j(x)+\frac{1}{2} K^j \|v\|^2,
\end{equation*}
which holds for each $i\in 	\langle p \rangle$, and $j\in	\langle m \rangle$. Therefore,
\begin{align}
	&\max_{i\in  	\langle p \rangle}h_j(x+v,\omega_i) \leq 	\digamma^*_j(x,v) + \digamma_j(x)+\frac{1}{2} K^j \|v\|^2, \notag 
	\\ \label{10}
	\text{i.e.,} ~&\digamma_j(x+v)\leq \digamma^*_j(x,v) + \digamma_j(x)+\frac{1}{2}K^j \|v\|^2.
\end{align}
The following  two results related to continuous approximation are true, for all $j \in \langle m \rangle$ \citet{bazaraa1982algorithm}:
\begin{enumerate}[($i$)]
	\item 	$\digamma^*_{j}(x,\lambda v)\leq \lambda \digamma^*_{j}(x,v),~\forall~\lambda \in [0,1]$.
	\item\label{m} $\digamma_j(x+\lambda v)\leq \digamma_j(x) +\lambda \digamma^*_j(x,v)+\frac{1}{2}\lambda^2 K^j \|v\|^2.$
\end{enumerate}
With the help of $(\ref{m})$ and Eqn. (\ref{10}), we can write
\begin{equation}\label{11}
	\digamma_j(x+\alpha v)\leq \alpha \digamma^*_j(x,v) + \digamma_j(x)+\frac{1}{2}\alpha^2K^j \|v\|^2,
\end{equation}
where $\alpha>0$ is sufficiently small. Since $x$ is not a critical point and $v\not=0,$ by Eqn. (\ref{11***}), we obtain 
\begin{center}
	$\digamma^*_j(x,v)< 0. $
\end{center}
For any $\beta \in (0,1),$ we  get
\begin{equation}\label{12}
	\digamma^*_j(x,v) < \beta  \digamma^*_j(x,v)
\end{equation}
As $\alpha$ is sufficiently small, the third term of the right hand side of the inequality (\ref{11}) tends to zero. Consequently, from inequalities (\ref{11}) and (\ref{12}), we get
\begin{equation}\label{13}
	\digamma_j(x+\alpha v)\leq  \digamma_j(x)+  \alpha \beta \digamma^*_j(x,v).
\end{equation}
By combining Eqn. (\ref{g}) and Eqn. (\ref{13}), we get
\begin{equation}\label{14*}
	\digamma_j(x+\alpha v)\leq  \digamma_j(x)+  \alpha \beta h(x,v).
\end{equation} 
Eqn. (\ref{14*}) denotes the step size rule for the conjugate gradient algorithm applied to $OWC_{P(U^F)}$.
\end{proof}

\par We have identified a descent direction and employed an inexact line search technique to determine the step length.~Subsequently, we formulate the conjugate gradient algorithm for $OWC_{P(U^F)}$ as follows.
\begin{alg1}\label{algo1}(Nonlinear conjugate gradient algorithm for $OWC_{P(U^F)}$)
\begin{enumerate}[{Step} 1]
	\item Choose $\epsilon>0$, $\beta \in (0,1) $   and $x^0\in \mathbb{R}^n$. Set $k :=0.$
	\item\label{step2} Solve $P(x^k)$ and find $s^k$ and $T({x^k})$.
	\item\label{step3} If  $|T(x^k)|<\epsilon$ or $\|s(x^k)\|<\epsilon$, then stop. Otherwise proceed to Step \ref{step4}.
	\item\label{step4} Define 
	\begin{align}\label{d}
		v^k =
		\begin{cases}
			s(x^k), & \text{if } k=0,\\
			s(x^k)+\gamma_kv^{k-1}, & \text{if }k\geq 1,
		\end{cases}
	\end{align}
	where $\gamma_k$ is an algorithmic parameter.
	\item \label{step5} Choose $\alpha_k$ as the largest $\alpha \in \{ \frac{1}{2^r} : r=1,2,3,\ldots\}$ satisfying Eqn.~(\ref{13}) and Eqn.~(\ref{14*}).
	\item Define $x^{k+1}:= x^k + \alpha_{k} v^k$, update $k:=k+1$ and go to Step \ref{step2}.
\end{enumerate}
\end{alg1}
\subsection{Justification of Algorithm \ref{algo1}}
In this section, we discuss the well-definedness of Algorithm \ref{algo1}, depending on Step \ref{step2}, Step \ref{step3}, Step \ref{step4}, and Step \ref{step5}. The key to the algorithm's well-definedness lies in the following steps:
\par In Step \ref{step2}, the algorithm requires computing a minimizer of the function $v \mapsto h(x^k, v) + \frac{1}{2}\|v\|^2$. Due to the strong convexity of this function, a unique minimizer exists, ensuring the existence of $s^k$. Thus, Step \ref{step2} is well-defined. Once $s^k$ is determined, $T(x^k)$ can be computed, affirming the well-definedness of Step \ref{step3}.
\par It is important to note, as per Theorem \ref{Thm 3.2}, that the stopping criterion $|T(x^k)| < \epsilon$ can be effectively substituted with $\|s^k\| < \epsilon$. This substitution is significant because if the algorithm does not proceed to Step \ref{step4} in iteration $k$ that is, if it stops at Step \ref{step3}, then Theorem \ref{Thm 3.2} implies that $x^k$ is an approximate critical point of $\digamma$. 
\par If Step \ref{step4} is reached at iteration $k$, to validation of Step \ref{step4} we have to choose $\gamma_k$ in such a way that the direction $v^k$ is a descent direction for $\digamma$ at $x^k.$ In order to validate the Step \ref{step4}, we presents the following lemma:
\begin{lem} Assume that in the Algorithm \ref{algo1}, the sequence $\gamma_k$ is defined so that it has the following property:\\
\begin{align}\label{p1}
	\gamma_k\in
	\begin{cases}
		[0,\infty), & ~\text{if }~h(x^k,v^{k-1})\leq0,\\
		[0,-a), & ~\text{if }~h(x^k,v^{k-1})>0,
	\end{cases}
\end{align}
or
\begin{align}\label{p2}
	\gamma_k\in
	\begin{cases}
		[0,\infty), & ~\text{if }~h(x^k,v^{k-1})\leq 0,\\
		[0,-\nu a], & ~\text{if }~h(x^k,v^{k-1})>0,
	\end{cases}
\end{align} for some $\nu\in[0,1),$ where $a=\frac{h(x^k,s(x^k))}{h(x^k,v^{k-1})}.$ If property given in Eqn. (\ref{p1}) holds, then $v^k$ is a descent direction for $\digamma.$ If property given in Eqn. (\ref{p2}) holds, then $v^k$ satisfies the sufficient descent condition with $b=1-\nu,$ i.e.,
\begin{equation}\label{sd}
	h(x^k,v^k)\leq b h(x^k,s(x^k)).
\end{equation}
\end{lem}
\begin{proof} Assume that $x^k$ is not a critical point, then 
\begin{align*}
	h(x^k,v^k)=\displaystyle \max_{	\langle m \rangle}\max_{i\in 	\langle p \rangle}  \{h_j(x^k,\omega_i)+\nabla h_j(x^k,\omega_i)^Tv^k-\digamma_j(x^k)\}.
\end{align*}
By Eqn. (\ref{d}), for $k=0,$ $v^k=s(x^0)$ hence the result. For $k\geq1,$ 
\begin{align*}
	h(x^k,d^k)&=\displaystyle \max_{		j \in \langle m \rangle}\max_{i\in 	\langle p \rangle}  \{h_j(x^k,\omega_i)+\nabla h_j(x^k,\omega_i)^T(s(x^k)+\gamma_kv^{k-1})-\digamma_j(x^k)\}\nonumber\\
	&\leq\displaystyle \max_{		j \in \langle m \rangle}\max_{i\in 	\langle p \rangle}  \{h_j(x^k,\omega_i)+\nabla h_j(x^k,\omega_i)^Ts(x^k)+\gamma_k\nabla h_j(x^k,\omega_i)^Tv^{k-1}-\digamma_j(x^k)\}\nonumber\\
	&\leq\displaystyle \max_{		j \in \langle m \rangle}\max_{i\in 	\langle p \rangle}  \{h_j(x^k,\omega_i)+\nabla h_j(x^k,\omega_i)^Ts(x^k)-\digamma_j(x^k)\}+\max_{		j \in \langle m \rangle}\max_{i\in 	\langle p \rangle} \gamma_k \nabla h_j(x^k,\omega_i)^Tv^{k-1}\nonumber\\
	&= h(x^k,s(x^k))+\gamma_k\max_{		j \in \langle m \rangle}\max_{i\in 	\langle p \rangle} \nabla h_j(x^k,\omega_i)^Tv^{k-1} \nonumber\\
	&\leq h(x^k,s(x^k))+\gamma_k\max_{		j \in \langle m \rangle}\max_{i\in 	\langle p \rangle}\{h_j(x^k,\omega_i)+  \nabla h_j(x^k,\omega_i)^Tv^{k-1}-\digamma_{j}(x)\}\nonumber\\
	&= h(x^k,s(x^k))+\gamma_kh(x^k,v^{k-1}).\nonumber
\end{align*}
Therefore,
\begin{equation}\label{e}
	h(x^k,d^k)\leq h(x^k,s(x^k))+\gamma_kh(x^k,v^{k-1}).
\end{equation}
Now, by property given in Eqn. (\ref{p1}) and Eqn. (\ref{p2}) if $h(x^k,v^{k-1})\leq0,$ then $\gamma_k\geq0$, and also $\gamma_kh(x^k,v^{k-1})\leq0.$ Then, by Eqn. (\ref{e}), $h(x^k,v^k)\leq h(x^k,s(x^k))$. Since $s(x^k)$ is a descent direction, according to Lemma \ref{lem1}, $h(x^k, s(x^k)) < 0$, and consequently, $v^k$ is also a descent direction. Also, if $h(x^k,v^{k-1})>0,$ then by property given in Eqn. (\ref{p2}), $\gamma_k\in[0,-\nu a],$ where $a=\frac{h(x^k,s(x^k))}{h(x^k,v^{k-1})}$ and $\nu\in[0,1)$. Now, if we take $\gamma_k=-\nu a=\frac{-\nu h(x^k,s(x^k))}{h(x^k,v^{k-1})},$ then by Eqn. (\ref{e}), we get 
\begin{align*}
	h(x^k,v^k)&\leq h(x^k,s(x^k))-\nu h(x^k,s(x^k))\\
	&\leq (1-\nu) h(x^k,s(x^k)).
\end{align*}
Thus, $v^k$ satisfies sufficient descent condition $h(x^k,v^k) \leq bh(x^k,s(x^k))$, with $b=1-\nu.$ Now, if $h(x^k,v^{k-1})>0,$ then by property given in Eqn. (\ref{p1}), $\gamma_k\in[0,-a).$ Therefore,
\begin{align*}\label{e}
	h(x^k,v^k)&\leq h(x^k,s(x^k))-ah(x^k,v^{k-1})\nonumber\\
	&=h(x^k,s(x^k))-\frac{h(x^k,s(x^k))}{h(x^k,v^{k-1})}h(x^k,v^{k-1})\nonumber\\
	&=0.
\end{align*}
Since $x^k$ is not a critical point and $h(x^k,v^k)\leq0,$  $h(x^k,v^k)$ can not be zero. Therefore, $h(x^k,v^k)<0$ implies that $v^k$ is a descent direction. Hence the Lemma. 
\end{proof}
With the help of above Lemma \ref{lem1}, it is clear that if $\gamma_k$ satisfies the property given in Eqn. (\ref{p1}) or Eqn. (\ref{p2}), then $v^k$ (defined in Eqn. \ref{d}) will be the the descent direction for $\digamma(x)$ at $x^k$ and hence Step \ref{step4} is well defined.\\
In Step \ref{step5},  we have two choices of $v^k,$ if $v^k=s(x^k),$ then in the direction $v^k,$ we choose $\alpha_k$ as the largest $\alpha\in  \{ \frac{1}{2^r} : r=1,2,3,\ldots\}$ which satisfies Eqn. (\ref{13}). If $v^k=s(x^k)+\gamma_kv^{k-1}$, and it is also proved that $h(x^k,v^k)<0$, in that case, $\alpha_k$ is chosen as the largest $\alpha\in  \{ \frac{1}{2^r} : r=1,2,3,\ldots\}$ that satisfies  \begin{equation}\label{linesearchwithd}
\digamma_j(x^k+\alpha_k v^k)\leq  \digamma_j(x^k)+  \alpha \beta h(x^k,v^k).
\end{equation}
Considering Eqn. (\ref{13}) and Eqn. (\ref{linesearchwithd}), it is evident that the objective function values consistently decreases in the component-wise partial order. This is due to the fact that, according to Theorem \ref{Thm 3.2}, $x^k$ is a noncritical point of $\digamma$, implying $T(x^k)<0$. Consequently, in Step \ref{step4}, it follows that $\digamma(x_{k+1})<\digamma(x_{k})$, ensuring the well-defined nature of Step \ref{step4}.
Utilizing $v^k$, $\alpha_k$, and the current iteration point $x^k$, we compute $x^{k+1}$ in Step \ref{step5}. Subsequently, we proceed to Step \ref{step2}, iterating through this process until the stopping criteria in Step \ref{step3} are met.\\ 
\par Before showing the convergence of the Algorithm \ref{algo1}, we discuss the another choice of $\gamma_k$ related to the $OWC_{P(U^F)}$ extension of FR, CD, DY, PRP, and HS for multiobjective optimization problems. For all the choices of $\gamma_k$, we will show that they all satisfy descent property and satisfies Armijo-type line search rule defined in Eqn. (\ref{14*}). 
\subsubsection{Fletcher-Reeves extension for $OWC_{P(U^F)}$}
The choice of $\gamma_k$ for $OWC_{P(U^F)}$ by Fletcher-Reeves can be extended as 
$\gamma^{FR}_k=\frac{h(x^k,s(x^k))}{h(x^{k-1},v(x^{k-1}))}.$
\begin{thm}
Let $\{x^k\}$ be a sequence generated by the Algorithm \ref{algo1} under consideration that $ 0\leq\gamma_k=\mu\gamma^{FR}_k$, where $0\leq\mu<1$ and suppose that $\alpha_{k}$ satisfies the Armijo-type inexact line search. Then, $v^k$ satisfies the sufficient descent condition (\ref{sd}) with $h(x^k,v^{k-1})<0$ and $b=1.$\end{thm}
\begin{proof}
Since $0\leq\gamma_k=\mu\gamma^{FR}_k$ and $\mu\geq0,$ $\gamma^{FR}_k\geq0$ and hence $\gamma^{FR}_k$ is well defined. 
By definition of $h(x,v)$ and Armijo-type inexact line search (\ref{14*}) we have 
\begin{align*}
	h(x^k,v^k)&=\displaystyle \max_{		j \in \langle m \rangle}\max_{i\in 	\langle p \rangle}  \{h_j(x^k,\omega_i)+\nabla h_j(x^k,\omega_i)^Tv^k-\digamma_j(x^k)\}\\
	&=\displaystyle \max_{		j \in \langle m \rangle}\max_{i\in 	\langle p \rangle}  \{h_j(x^k,\omega_i)+\nabla h_j(x^k,\omega_i)^T(s(x^k)+\gamma_kv^{k-1})-\digamma_j(x^k)\}\\
	&=\displaystyle \max_{		j \in \langle m \rangle}\max_{i\in 	\langle p \rangle}  \{h_j(x^k,\omega_i)+\nabla h_j(x^k,\omega_i)^Ts(x^k)+\gamma_kh_j(x^k,\omega_i)^Tv^{k-1}-\digamma_j(x^k)\}\\
	&\leq\displaystyle \max_{		j \in \langle m \rangle}\max_{i\in 	\langle p \rangle}  \{h_j(x^k,\omega_i)+\nabla h_j(x^k,\omega_i)^Ts(x^k)-\digamma_j(x^k)\}+\displaystyle \max_{		j \in \langle m \rangle}\max_{i\in 	\langle p \rangle}\gamma_kh_j(x^k,\omega_i)^Tv^{k-1}\\
	&=h(x^k,s(x^k))+\gamma_k\displaystyle \max_{		j \in \langle m \rangle}\max_{i\in 	\langle p \rangle}\nabla h_j(x^k,\omega_i)^Tv^{k-1}\\
	&\leq h(x^k,s(x^k))+\gamma_k\displaystyle \max_{		j \in \langle m \rangle}\max_{i\in 	\langle p \rangle}\{h_j(x^k,\omega_i)+\nabla h_j(x^k,\omega_i)^Tv^{k-1}-\digamma_j(x)\}\\
	&= h(x^k,s(x^k))+\gamma_kh(x^k,v^{k-1})\\
	&= h(x^k,s(x^k))+\mu\gamma_kh(x^k,v^{k-1})\\
	&\leq h(x^k,s(x^k))+\mu\gamma^{FR}_kh(x^k,v^{k-1})\\
	&= h(x^k,s(x^k))+\mu\frac{h(x^k,s(x^k))}{h(x^{k-1},v(x^{k-1}))}h(x^k,v^{k-1}),
\end{align*}
since $s(x^k)$ and $v(x^{k-1})$ are the descent directions of $\digamma$ at $x^k$ and $x^{k-1}$, respectively, and $h(x^k,v^{k-1})<0$, therefore, $h(x^k,v^k)\leq h(x^k,s(x^k))$ provided $b=1$, i.e., the sufficient descent condition is satisfied, for all $k\geq0.$ Hence, $v^{k}$ is a descent direction.
\end{proof}
\subsubsection{Conjugate descent extension for $OWC_{P(U^F)}$}
The choice of $\gamma_k$ for $OWC_{P(U^F)}$ of conjugate descent can be extended as 
$\gamma^{CD}_k=\frac{h(x^k,s(x^k))}{h(x^{k-1},v^{k-1})}.$\\
Now, in the next theorem we prove that the Algorithm \ref{algo1} generates a descent direction if $ 0\leq\gamma_k\leq\gamma^{CD}_k$ and the $\alpha_{k}$ satisfies the Armijo-type inexact line search.
\begin{thm}
Let $\{x^k\}$ be a sequence generated by the Algorithm \ref{algo1} under consideration that $ 0\leq\gamma_k\leq\gamma^{CD}_k$, and suppose that $\alpha_{k}$ satisfies the Armijo-type inexact line search. Then, $v^k$ satisfies the sufficient descent condition (\ref{sd}) with $b=1-\mu.$\end{thm}
\begin{proof}
Proof of the theorem can be given by induction. At $k=0$, $v^0=s(x^0),$ i.e., $v^0$ is a descent direction, i.e.,  $h(x^0,s(x^0))<0.$ If $0<\mu<1$, then Eqn. (\ref{sd}) with $b=1-\mu$ holds. Now, for $k\geq1$ assume that \begin{equation}
	h(x^{k-1},v^{k-1})\leq(1-\mu)h(x^{k-1},v(x^{k-1}))<0.
\end{equation}
Then, $\gamma^{CD}_k>0$ and $\gamma_k$ is well defined. Also, by definition of $h(x,v)$ and Armijo-type inexact line search technique given in (\ref{14*}), we have 
\begin{align*}
	h(x^k,v^k)&=\displaystyle \max_{		j \in \langle m \rangle}\max_{i\in 	\langle p \rangle}  \{h_j(x^k,\omega_i)+\nabla h_j(x^k,\omega_i)^Tv^k-\digamma_j(x^k)\}\\
	&=\displaystyle \max_{		j \in \langle m \rangle}\max_{i\in 	\langle p \rangle}  \{h_j(x^k,\omega_i)+\nabla h_j(x^k,\omega_i)^T(s(x^k)+\gamma_kv^{k-1})-\digamma_j(x^k)\}\\
	&=\displaystyle \max_{		j \in \langle m \rangle}\max_{i\in 	\langle p \rangle}  \{h_j(x^k,\omega_i)+\nabla h_j(x^k,\omega_i)^Ts(x^k)+\gamma_kh_j(x^k,\omega_i)^Tv^{k-1}-\digamma_j(x^k)\}\\
	&\leq\displaystyle \max_{		j \in \langle m \rangle}\max_{i\in 	\langle p \rangle}  \{h_j(x^k,\omega_i)+\nabla h_j(x^k,\omega_i)^Ts(x^k)-\digamma_j(x^k)\}+\displaystyle \max_{		j \in \langle m \rangle}\max_{i\in 	\langle p \rangle}\gamma_kh_j(x^k,\omega_i)^Tv^{k-1}\\
	&=h(x^k,s(x^k))+\gamma_k\displaystyle \max_{		j \in \langle m \rangle}\max_{i\in 	\langle p \rangle}\nabla h_j(x^k,\omega_i)^Tv^{k-1}\\
	&\leq h(x^k,s(x^k))+\gamma_k\displaystyle \max_{		j \in \langle m \rangle}\max_{i\in 	\langle p \rangle}\{h_j(x^k,\omega_i)+\nabla h_j(x^k,\omega_i)^Tv^{k-1}-\digamma_j(x)\}\\
	&= h(x^k,s(x^k))+\gamma_kh(x^k,v^{k-1})\\
	&\leq h(x^k,s(x^k))-\mu\gamma_kh(x^k,v^{k-1})\\
	&\leq h(x^k,s(x^k))-\mu\gamma^{CD}_kh(x^k,v^{k-1})\\
	&= (1-\mu)h(x^k,s(x^k)).
\end{align*}
Therefore, $h(x^{k},v^{k})\leq(1-\mu)h(x^k,s(x^k)),$ i.e., sufficient descent condition is satisfy, for all $k\geq 0.$ Hence, $v^{k}$ is a descent direction.
\end{proof}
\subsubsection{Dai-Yuan extension for $OWC_{P(U^F)}$}
The choice of $\gamma_k$ for $OWC_{P(U^F)}$ by Dai-Yuan can be extended as 
$$\gamma^{DY}_k=\frac{-h(x^k,s(x^k))}{f(x^{k},v^{k-1})-h(x^{k-1},v^{k-1})}.$$
\begin{thm}
Let $\{x^k\}$ be a sequence generated by the Algorithm \ref{algo1} under consideration that $ 0\leq\gamma_k\leq\gamma^{DY}_k$, and suppose that $\alpha_{k}$ satisfies the Armijo-type inexact line search. Then, $v^k$ satisfies the sufficient descent condition (\ref{sd}) with $b=\frac{1}{1+\mu}.$
\end{thm}
\begin{proof}
Proof of the theorem can be given by induction. At $k=0$, $v^0=s(x^0)$, i.e., $v^0$ is a descent direction, i.e., $h(x^0,v^0)<0$. If $0<\mu<1$, then Eqn. (\ref{sd}) with $b=\frac{1}{1+\mu}$ holds. Now, for some $k\geq 1$ assume that \begin{equation}
	h(x^{k-1},v^{k-1})\leq\frac{1}{1+\mu}h(x^{k-1},v(x^{k-1}))<0. 
\end{equation}
By Armijo-type inexact line search (\ref{14*}), we obtain 
$$h(x^k,v^{k-1})\geq \mu f(x^{k-1},v^{k-1})>h(x^{k-1},v^{k-1}),$$
because $\mu<1$ and $h(x^{k-1},v^{k-1})<0.$ Therefore, $\gamma_k\leq\gamma^{DY}_k>0$ and $\gamma_k$ is well defined. 
Also, by definition of $h$, $v^k$ and positiveness of $\gamma_k$ we have 
\begin{align*}
	h(x^k,v^k)&=\displaystyle \max_{		j \in \langle m \rangle}\max_{i\in 	\langle p \rangle}  \{h_j(x^k,\omega_i)+\nabla h_j(x^k,\omega_i)^Tv^k-\digamma_j(x^k)\}\\
	&=\displaystyle \max_{		j \in \langle m \rangle}\max_{i\in 	\langle p \rangle}  \{h_j(x^k,\omega_i)+\nabla h_j(x^k,\omega_i)^T(s(x^k)+\gamma_kv^{k-1})-\digamma_j(x^k)\}\\
	&=\displaystyle \max_{		j \in \langle m \rangle}\max_{i\in 	\langle p \rangle}  \{h_j(x^k,\omega_i)+\nabla h_j(x^k,\omega_i)^Ts(x^k)+\gamma_kh_j(x^k,\omega_i)^Tv^{k-1}-\digamma_j(x^k)\}\\
	&\leq\displaystyle \max_{		j \in \langle m \rangle}\max_{i\in 	\langle p \rangle}  \{h_j(x^k,\omega_i)+\nabla h_j(x^k,\omega_i)^Ts(x^k)-\digamma_j(x^k)\}+\displaystyle \max_{		j \in \langle m \rangle}\max_{i\in 	\langle p \rangle}\gamma_kh_j(x^k,\omega_i)^Tv^{k-1}\\
	&=h(x^k,s(x^k))+\gamma_k\displaystyle \max_{		j \in \langle m \rangle}\max_{i\in 	\langle p \rangle}\nabla h_j(x^k,\omega_i)^Tv^{k-1}\\
	&\leq h(x^k,s(x^k))+\gamma_k\displaystyle \max_{		j \in \langle m \rangle}\max_{i\in 	\langle p \rangle}\{h_j(x^k,\omega_i)+\nabla h_j(x^k,\omega_i)^Tv^{k-1}-\digamma_j(x)\}\\
	&= h(x^k,s(x^k))+\gamma_kh(x^k,v^{k-1}).
\end{align*}
Therefore,
$h(x^k,v^k)\leq h(x^k,s(x^k))+\gamma_kh(x^k,v^{k-1}).$ Now if $h(x^k,v^{k-1})\leq 0,$ then we get the result. If we assume $h(x^k,v^{k-1})>0$. Then, \begin{align}\label{dy}
	h(x^k,v^k)\leq h(x^k,s(x^k))+\gamma_kh(x^k,v^{k-1})\leq h(x^k,s(x^k))+\gamma^{DY}_kh(x^k,v^{k-1}).
\end{align}
Take $w_k=\frac{h(x^k,v^{k-1})}{h(x^{k-1},v^{k-1})}$ and consider $w_k\in[-\mu,\mu].$ Using the Eqn. (\ref{dy}), we obtain 
\begin{align*}
	h(x^k,v^k)\leq \frac{1}{1-w_k}h(x^k,s(x^k)\leq  h(x^k,s(x^k))\leq\frac{1}{1+\mu}h(x^k,s(x^k)),
\end{align*}
since $h(x^k,s(x^k))<0$ and $\frac{1}{1+\mu}>0$ then by above inequality $h(x^k,v^k)<0$ hence the proof. 
\end{proof}
\subsubsection{Polak-Ribi$\grave{e}$re-Polyak extension for $OWC_{P(U^F)}$}
The choice of $\gamma_k$ for $OWC_{P(U^F)}$ by Polak-Ribiere-Polak can be extended as 
$$\gamma^{PRP}_k=\frac{-h(x^k,s(x^k))+h(x^{k-1},v(x^{k}))}{-h(x^{k-1},v(x^{k-1}))}.$$
\begin{thm}
Let $\{x^k\}$ be a sequence generated by the Algorithm \ref{algo1} under consideration that $ \gamma_k=\max\{\gamma^{PRP}_k,0\}$, and suppose that $\alpha_{k}$ satisfies the Armijo-type inexact line search. Then, $v^k$ satisfies the sufficient descent condition (\ref{sd}) provided $h(x^k,s(x^k))<h(x^{k-1},v(x^{k})).$
\end{thm}
\begin{proof}
By definition of $h$, $v^k$ and positiveness of $\gamma_k$ we have 
\begin{align*}
	h(x^k,v^k)&=\displaystyle \max_{		j \in \langle m \rangle}\max_{i\in 	\langle p \rangle}  \{h_j(x^k,\omega_i)+\nabla h_j(x^k,\omega_i)^Tv^k-\digamma_j(x^k)\}\\
	&=\displaystyle \max_{		j \in \langle m \rangle}\max_{i\in 	\langle p \rangle}  \{h_j(x^k,\omega_i)+\nabla h_j(x^k,\omega_i)^T(s(x^k)+\gamma_kv^{k-1})-\digamma_j(x^k)\}\\
	&=\displaystyle \max_{		j \in \langle m \rangle}\max_{i\in 	\langle p \rangle}  \{h_j(x^k,\omega_i)+\nabla h_j(x^k,\omega_i)^Ts(x^k)+\gamma_kh_j(x^k,\omega_i)^Tv^{k-1}-\digamma_j(x^k)\}\\
	&\leq\displaystyle \max_{		j \in \langle m \rangle}\max_{i\in 	\langle p \rangle}  \{h_j(x^k,\omega_i)+\nabla h_j(x^k,\omega_i)^Ts(x^k)-\digamma_j(x^k)\}+\displaystyle \max_{		j \in \langle m \rangle}\max_{i\in 	\langle p \rangle}\gamma_kh_j(x^k,\omega_i)^Tv^{k-1}\\
	&=h(x^k,s(x^k))+\gamma_k\displaystyle \max_{		j \in \langle m \rangle}\max_{i\in 	\langle p \rangle}\nabla h_j(x^k,\omega_i)^Tv^{k-1}\\
	&\leq h(x^k,s(x^k))+\gamma_k\displaystyle \max_{		j \in \langle m \rangle}\max_{i\in 	\langle p \rangle}\{h_j(x^k,\omega_i)+\nabla h_j(x^k,\omega_i)^Tv^{k-1}-\digamma_j(x)\}\\
	&= h(x^k,s(x^k))+\gamma_kh(x^k,v^{k-1}).
\end{align*}
Therefore,
\begin{equation}\label{prp}
	h(x^k,v^k)\leq h(x^k,s(x^k))+\gamma_kh(x^k,v^{k-1}).
\end{equation} 
Now, if $h(x^k,v^{k-1})\leq 0,$ then we obtain the result. If we assume $h(x^k,v^{k-1})>0$. Then, by $\gamma_k=\max\{\gamma^{PRP}_k,0\},$ two cases may be arise:
$$Case~1:~~~\text{If}~ \gamma_k^{PRP}\leq 0, ~\text{then}~\gamma_k=0$$
$$Case~2:~~~\text{If}~ \gamma_k^{PRP}>0, ~\text{then}~\gamma_k=\gamma^{PRP}.$$
If we consider $Case~1,$ then $h(x^k,v^k)\leq h(x^k,s(x^k)),$ i.e., sufficient decent condition is satisfied and $v^k$ is a descent direction. Also, if we consider $Case~2,$ then by Eqn. (\ref{prp}), we have 
\begin{align*}
	h(x^k,v^k)&\leq h(x^k,s(x^k))+\gamma^{PRP}_kh(x^k,v^{k-1})\\
	&=h(x^k,s(x^k))+\bigg(\frac{-h(x^k,s(x^k))+h(x^{k-1},v(x^{k}))}{-h(x^{k-1},v(x^{k-1}))}\bigg)h(x^k,v^{k-1})\\
	&=h(x^k,s(x^k))+\bigg(\frac{h(x^k,s(x^k))}{h(x^{k-1},v(x^{k-1}))}-\frac {h(x^{k-1},v(x^{k}))}{h(x^{k-1},v(x^{k-1}))}\bigg)h(x^k,v^{k-1})
\end{align*}
Since $h(x^{k-1},s(x^k))>h(x^{k},s(x^k))$ and $h(x^k,v^{k-1})>0,$ the second term of the right hand side of the above inequality will be negative. Thus, $h(x^k,v^k)\leq h(x^k,s(x^k)),$ i.e., sufficient decent condition is satisfied and $v^k$ is a descent direction.
\end{proof}
\subsubsection{Hestenes-Stiefel extension for $OWC_{P(U^F)}$}
The choice of $\gamma_k$ for $OWC_{P(U^F)}$ by Hestenes-Stiefel can be extended as 
$$\gamma^{HS}_k=\frac{-h(x^k,s(x^k))+h(x^{k-1},v(x^{k}))}{h(x^k,v^{k-1}))-h(x^{k-1},v^{k-1}))}.$$
\begin{thm}
Let $\{x^k\}$ be a sequence generated by the Algorithm \ref{algo1} under consideration that $ \gamma_k=\max\{\gamma^{HS}_k,0\}$, and suppose that $\alpha_{k}$ satisfies the Armijo-type inexact line search. Then, $v^k$ satisfies the sufficient descent condition (\ref{sd}) provided $h(x^k,s(x^k))<h(x^{k-1},v(x^{k}))$ and $h(x^k,v^{k-1}))<h(x^{k-1},v^{k-1})).$
\end{thm}
\begin{proof}
By definition of $h$, $v^k$ and positiveness of $\gamma_k$ we have 
\begin{align*}
	h(x^k,v^k)&=\displaystyle \max_{		j \in \langle m \rangle}\max_{i\in 	\langle p \rangle}  \{h_j(x^k,\omega_i)+\nabla h_j(x^k,\omega_i)^Tv^k-\digamma_j(x^k)\}\\
	&=\displaystyle \max_{		j \in \langle m \rangle}\max_{i\in 	\langle p \rangle}  \{h_j(x^k,\omega_i)+\nabla h_j(x^k,\omega_i)^T(s(x^k)+\gamma_kv^{k-1})-\digamma_j(x^k)\}\\
	&=\displaystyle \max_{		j \in \langle m \rangle}\max_{i\in 	\langle p \rangle}  \{h_j(x^k,\omega_i)+\nabla h_j(x^k,\omega_i)^Ts(x^k)+\gamma_kh_j(x^k,\omega_i)^Tv^{k-1}-\digamma_j(x^k)\}\\
	&\leq\displaystyle \max_{		j \in \langle m \rangle}\max_{i\in 	\langle p \rangle}  \{h_j(x^k,\omega_i)+\nabla h_j(x^k,\omega_i)^Ts(x^k)-\digamma_j(x^k)\}+\displaystyle \max_{		j \in \langle m \rangle}\max_{i\in 	\langle p \rangle}\gamma_kh_j(x^k,\omega_i)^Tv^{k-1}\\
	&=h(x^k,s(x^k))+\gamma_k\displaystyle \max_{		j \in \langle m \rangle}\max_{i\in 	\langle p \rangle}\nabla h_j(x^k,\omega_i)^Tv^{k-1}\\
	&\leq h(x^k,s(x^k))+\gamma_k\displaystyle \max_{		j \in \langle m \rangle}\max_{i\in 	\langle p \rangle}\{h_j(x^k,\omega_i)+\nabla h_j(x^k,\omega_i)^Tv^{k-1}-\digamma_j(x)\}\\
	&= h(x^k,s(x^k))+\gamma_kh(x^k,v^{k-1}).
\end{align*}
Therefore,
\begin{equation}\label{hs}
	h(x^k,v^k)\leq h(x^k,s(x^k))+\gamma_kh(x^k,v^{k-1}).
\end{equation} 
Now if $h(x^k,v^{k-1})\leq 0,$ then the result is obtained. If we assume $h(x^k,v^{k-1})>0$. Then, by $\gamma_k=\max\{\gamma^{HS}_k,0\}$ two cases may be arise:
$$Case~1:~~~\text{If}~ \gamma_k^{HS}\leq 0 ~\text{then}~\gamma_k=0$$
$$Case~2:~~~\text{If}~ \gamma_k^{HS}>0 ~\text{then}~\gamma_k=\gamma_k^{HS}.$$
If we consider $Case~1,$ then $h(x^k,v^k)\leq h(x^k,s(x^k)),$ i.e., sufficient decent condition is satisfied and $v^k$ is a descent direction. Also, if we consider $Case~2,$ then by Eqn. (\ref{hs}), we have 
\begin{align*}
	h(x^k,v^k)&\leq h(x^k,s(x^k))+\gamma^{HS}_kh(x^k,v^{k-1})\\
	&=h(x^k,s(x^k))+\bigg(\frac{-h(x^k,s(x^k))+h(x^{k-1},v(x^{k}))}{h(x^k,v^{k-1}))-h(x^{k-1},v^{k-1}))}\bigg)h(x^k,v^{k-1})\\
	&=h(x^k,s(x^k))+\bigg(\frac{h(x^k,s(x^k))}{h(x^{k-1},v(x^{k-1}))}-\frac {h(x^{k-1},v(x^{k}))}{h(x^{k-1},v(x^{k-1}))}\bigg)h(x^k,v^{k-1}).
\end{align*}

Since $h(x^{k-1},s(x^k)) > h(x^k,s(x^k)),$ $h(x^k,v^{k-1}) < h(x^{k-1},v^{k-1}),$ and $h(x^k,v^{k-1}) > 0,$ the second term on the right-hand side of the above inequality will be negative. Thus, $h(x^k,v^k) \leq h(x^k,s(x^k)),$ i.e., the sufficient descent condition is satisfied, and $v^k$ is a descent direction.
\end{proof}
\subsection{Convergence analysis of nonlinear conjugate gradient method (Algorithm \ref{algo1})}\label{ss3.3}
It is obvious that if Algorithm \ref{algo1} has finite iteration, then the last iterative point is approximately a critical point, and therefore it is an optimum for $\digamma$. So, it is relevant to consider the convergence analysis when Algorithm \ref{algo1} generates an infinite sequence. In view of this consideration, we assume that the $\{x^k\}$, $\{v^k\}$ and $\{\alpha_k\}$ are infinite sequences gererated by Algorithm \ref{algo1} for $OWC_{P(U^F)}$. We show that any accumulation point of $\{x^k\}$ is a critical point for $\digamma$.
\begin{thm}
Let $\{x^k\}$ be a sequence which is produced by Algorithm \ref{algo1}. If any accumulation point of $\{x^k\}$ is exists, then it will be the critical point for $\digamma$.
\end{thm}
\begin{proof}
Let  $x^*$ be an accumulation point  of the sequence $\{x^k\}$. The convergence of the $\{x^k\}$ depends on $v^k,$ where
\begin{align}\label{d}
	v^k =
	\begin{cases}
		s(x^k), & \text{if } k=0,\\
		s(x^k)+\gamma_kv^{k-1}, & \text{if }k\geq 1.
	\end{cases}
\end{align}
If we consider the case $v^k=s(x^k),$ where $s(x^k)$ is the solution of subproblem (\ref{eq31}) at $x^k.$ 
The solution and optimal value of the subproblem \eqref{eq31} at $x=x^*$ are given by
\begin{center}
	~~~~~~~~$s(x^*)	= \underset{v\in \mathbb{R}^n}{\mathrm{argmin}} \left\{  \displaystyle h(x^*,v) +\frac{1}{2} \|v\|^2\right\}$
	\text{and}~~$T(x^*)	= h(x^*,s(x^*))   +\frac{1}{2} \|s(x^*)\|^2$
\end{center}
respectively, where $h(x^*,v)=\displaystyle\max_{		j \in \langle m \rangle}\max_{i\in	\langle p \rangle}\{ h_j(x^*,\omega_i) + \nabla h_j(x^*,\omega_i)^Tv-\digamma_j(x^*) \}.$ We have to show $x^*$ is a critical point, by Theorem \ref{Thm 3.2}, that is, it suffices to show $T(x^*)=0$ or $v(x^*)=0$ or $h(x^*,s(x^*))=0,$ i.e., there is no descent direction at the point $x^*.$ By Algorithm \ref{algo1},  we know
$\digamma(x^k)$ is $\mathbb{R}_+^m$-decreasing, i.e., component-wise decreasing. Therefore, 
\begin{center}
	$\lim\limits_{k\rightarrow\infty}\digamma(x^k) = \digamma(x^*),$ as $x^k\to x^*$, because of continuity of $\digamma.$
\end{center}
Which implies
\begin{center}
	$\lim\limits_{k\rightarrow\infty}\|\digamma(x^k) - \digamma(x^*)\|=0.$
\end{center}
In the sense of component-wise we will get
\begin{center}
	$ \lim\limits_{k\rightarrow\infty}\digamma_{j}(x^k) = \digamma_{j}(x^*), ~\text{for all}~~j=1,2,\ldots,m.$
\end{center}
Now, by Algorithm \ref{algo1}, $Step~5,$ and Eqn. (\ref{14*}) we have

$$	\digamma_j(x^k+\alpha_k s^k)\leq \digamma_j(x^k)+\alpha_k \beta h(x^k,s^k)$$
and	\begin{equation}\label{15}
	\digamma_j(x^k)-\digamma_j(x^k+\alpha_k s^k) \geq -\alpha_k \beta h(x^k,s^k).
\end{equation}
As we know $\alpha_k$ and $\beta_k$ are positive, also $x^k$ is not a critical point and $s^k=s(x^k)$ is a descent direction that implies $h(x^k,s^k)<0$ and by Eqn. (\ref{15}), we get 
\begin{equation}\label{15}
	0\leq \digamma_j(x^k)-\digamma_j(x^k+\alpha_k s^k) \geq -\alpha_k \beta h(x^k,s^k)\geq 0.
\end{equation}
Therefore,
\begin{equation}\label{14}
	\displaystyle \lim\limits_{k\rightarrow\infty}\alpha_k h(x^k,s^k)=0, ~\text{for all}~ j \in\Lambda.
\end{equation}
Assuming that $\alpha_{k}\in(0,1],~\forall ~k\geq0$, then the following two cases may be arise:\\
$$Case~1:~~\limsup\limits_{k\rightarrow\infty} \alpha_{k}>0 $$  $$Case~2:~~ \limsup\limits_{k\rightarrow\infty} \alpha_{k}=0.$$
If we take the Case~1, then there exists a subsequence $\{x^{k_l}\}$ of $\{x^k\}$ converging to $x^*,$ and $\alpha^*>0$ such that $\lim\limits_{l\rightarrow\infty}\alpha_{k_l}=\alpha^*$. Then, by  Eqn. (\ref{14})
we get
\begin{equation*}
	0=\lim\limits_{l\rightarrow\infty}h(x^{k_l},s^{k_l})
	\leq \lim\limits_{l\rightarrow\infty} \bigg( h(x^{k_l},s^{k_l}) +\tfrac{1}{2} \|s(x^{k_l})\|^2\bigg),
\end{equation*}
and therefore
$$0\leq \lim\limits_{l\rightarrow\infty}T(x^{k_l}).$$
Since $T(x)\leq 0,$ for each $x$, we get
$$0 \leq T(x^*) \leq 0,$$
which implies $T(x^*)= 0$ and hence $x^*$ is a critical point.\\
Now, if we consider Case~2 and Theorem \ref{Thm 3.2}, $s(x)$ is bounded on any compact set. In other words, every term of the sequence ${s^k}$ lies in a compact set. Therefore, for each $k$, $s^k = s(x^k)$ lies in a bounded set. Consequently, ${s^k}$ is a bounded sequence on $C$. Due to the boundedness of $\{s^k\}$ there must exists a subsequence $\{s^{k_l}\}$ of $\{s^k\}$ such that $\lim\limits_{l\rightarrow\infty}s^{k_l}=s^*$ and $\lim\limits_{l\rightarrow\infty}\alpha_{k_l}=0.$\\
Note that we have
$$h(x^{k_l},s^{k_l})\leq T(x^{k_l})<0, ~\forall ~l\geq0.$$

So, as $l\to\infty$ we get
\begin{equation}\label{17}
	h(x^*,s^*)s^*\leq T(x^*)\leq 0.
\end{equation}
Since $\alpha_{k_l} \rightarrow 0$ for $l$ large enough, by Archimedian  property, $\alpha_{k_l}<\frac{1}{2^r}$, which means that the Armijo-like line search rule is not satisfied for $\alpha = \frac{1}{2^r},$ i.e.,
\begin{center}
	$	\digamma_j(x^{k_l}+\frac{1}{2^r} s^{k_l})\not\leq \digamma_j(x^{k_l})+ \frac{1}{2^r}\beta h(x^{k_l},s^{k_l}).$
\end{center}
So, for all $j\in\langle m \rangle$, there exists $j=j(k_l)\in \langle m \rangle$ such that
\begin{center}
	$	\digamma_j(x^{k_l}+\frac{1}{2^r} s^{k_l})\geq \digamma_j(x^{k_l})+ \frac{1}{2^r}\beta h(x^{k_l},s^{k_l}).$
\end{center}
Since  $\{j(k_l)\}_l \in \langle m \rangle$, there exist a subsequence $\{ k_{l_z}\}_z$
and an index $j_0$ such that\\ $j_0=j(k_{l_z}),~\forall~ z \geq1$  and
\begin{center}
	$	\digamma_j(x^{k_{l_z}}+\frac{1}{2^r} s^{k_{l_z}}) \geq \digamma_j(x^{k_{l_z}}) + \frac{1}{2^r}\beta \digamma^*_j({} x^{k_{l_z}},s^{k_{l_z}}). $
\end{center}
Taking limit $z \rightarrow \infty$ in the above inequality we obtain
\begin{center}
	$\digamma_j(x^*+\frac{1}{2^r} s^*) \geq \digamma_j(x^*) + \frac{1}{2^r}\beta h( x^*,s^*).$
\end{center}
Since this inequality  holds for any positive integer $r$ and for $j_0$ (depending on $r$), from Eqn. (\ref{13}), we obtain
\begin{center}
	$	\digamma_j(x+\alpha v)\leq  \digamma_j(x)+  \alpha \beta h(x,v).$
\end{center}
It follows that
\begin{center}
	$h( x^*,s^*)\not <0$.
\end{center}
So,
\begin{equation}\label{18}
	h( x^*,s^*)\geq 0.
\end{equation}
From Eqn. (\ref{17}) and Eqn. (\ref{18}), $T (x^*)=0$. Therefore, we  can conclude that $x^*$ is a critical point. Since we know $v^k=	s(x^k)+\gamma_kv^{k-1}$ is a descent direction for every choice of $\gamma_k.$ By simply replacing $s(x^k)$ with $v^k$, it is then possible to demonstrate the convergence of the NCGM in the case of $v^k$.  
\end{proof}
\section{Numerical results}\label{sec5}
In this section, we present numerical results related to the NCGM (Algorithm \ref{algo1}). We compare Algorithm \ref{algo1} with the weighted sum method based on the Pareto front, evaluating performance profiles in terms of the number of iterations, function evaluations, the $\Delta$ spread metric, and the hypervolume metric. Different scalarization approaches have been proposed in the literature, and their effectiveness can vary based on the nature of the objectives and constraints. Here are some commonly used scalarization methods: the weighted sum method, Tchebycheff approach, boundary intersection method, goal programming, and the $\epsilon-$constraint method. In most cases, it is worth noting that the weighted sum method often yields better results compared to other approaches. That is why we choose the weighted sum method for comparison with NCGM. However, the weighted sum method has some drawbacks that will be addressed by the NCGM.
We implemented Algorithm \ref{algo1} in Python. The search direction at the initial step is calculated by solving the subproblem $P(x^0)$ (i.e., $P(x)$ at $x^0$) using the $cvxpt$ solver. With the help of $P(x^0),$ $s(x^0)$ and $T(x^0)$ are calculated. For computing a step size, we used Armijo-type inexact line search techniques given in Eqn. \ref{13} and Eqn. \ref{14*}. In our computations, we consider $\epsilon=10^{-4}$ as tolerance or maximum $5000$ number of iterations is considered as stopping criteria. Thus, in our computations, we use $\|v^k\|<10^{-4}$ or $|T(x^k)|<10^{-4}$ or maximum $5000$ number of iterations as a stopping criteria.  Note that, if the stopping criterion is met at $x^0$, then $x^0$ is a critical point; otherwise, we proceed to the next step. Now we calculate $x^1$ with the help of $s(x^0)$ and step size $\alpha_0$ then we calculate $s(x^1)$ by solving $P(x^1)$. By using $v^0=s(x^0)$ and $s(x^1),$  we find $v^1=s(x^1)+\beta v^0$ and repeat this process until the stopping criterion is met. In the weighted sum method, we solve the following single-objective optimization problem
$$\min_{x\in\mathbb{R}^n} \left(a_1\digamma_1(x)+a_2\digamma_2(x)+\dots+a_m\digamma_m(x)\right),$$ where $a=(a_1,a_2,\dots,a_m)$ such that $a_i\geq0,$ using the technique developed in \citet{bazaraa1982algorithm} with initial approximation $x^0=\frac{1}{2}(lb+ub)$. For bi-objective optimization problems, we have considered weights $(1,0)$, $(0,1)$, and 98 random weights uniformly distributed in the square area of~$[0,1]\times[0,1]$ (i.e., any 98 random weights uniformly distributed in this area).~On the other hand, for three-objective optimization problems, we have considered four types of weights: $(1,0,0)$, $(0,1,0)$, $(0,0,1)$, and 97 random weights uniformly distributed in the cubic area of $[0,1]\times[0,1]\times[0,1]$ (i.e., any 97 random weights uniformly distributed in this area). 
\par
\begin{figure}
	\centering       
	\begin{subfigure}{1.2in}
		\includegraphics[width=1.1in,height=1.3in]
		{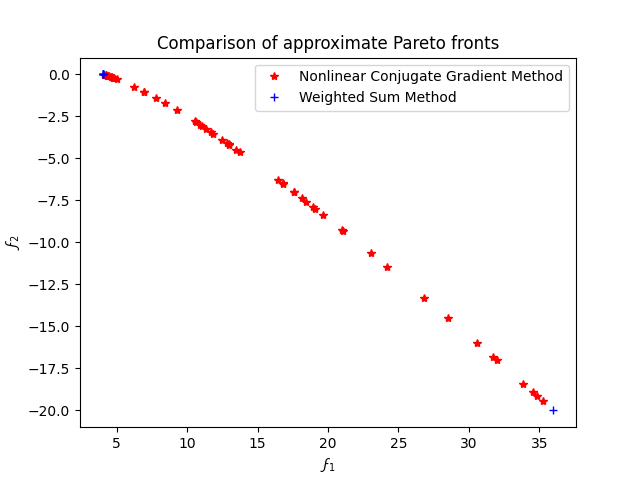}
		\centering    \caption*{($a_1$)}
	\end{subfigure}%
	\begin{subfigure}{1.2in}
		\includegraphics[width=1.1in,height=1.3in]
		{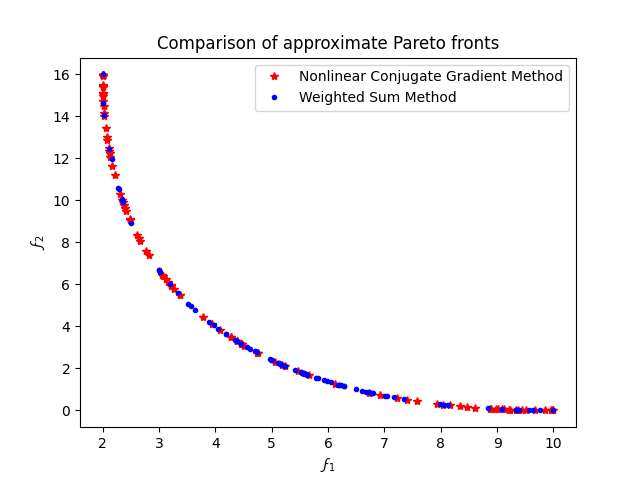}
		\centering     \caption*{($a_2$)}
	\end{subfigure}%
	\begin{subfigure}{1.2in}
		\includegraphics[width=1.1in,height=1.3in]
		{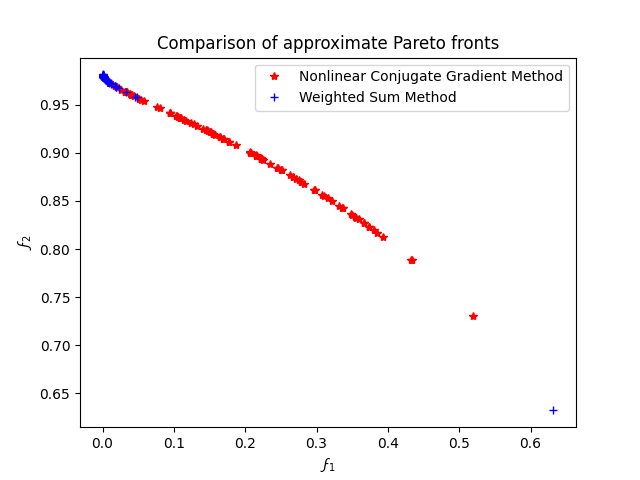}
		\centering   \caption*{($a_3$)}
	\end{subfigure}
	\begin{subfigure}{1.2in}
		\includegraphics[width=1.1in,height=1.3in]
		{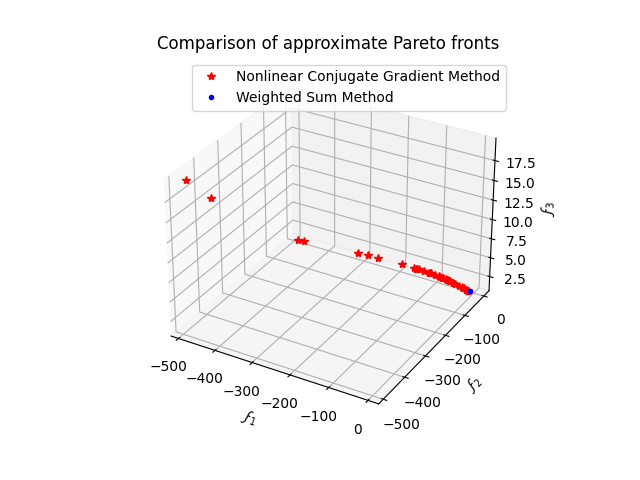}
		\centering   \caption*{($a_4$)}
	\end{subfigure}\\
	\begin{subfigure}{1.2in}
		\includegraphics[width=1.1in,height=1.3in]
		{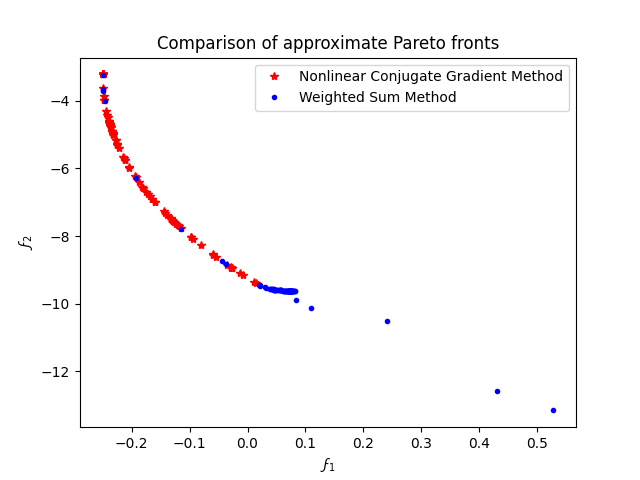}
		\centering   \caption*{($a_5$)}
	\end{subfigure}
	\begin{subfigure}{1.2in}
		\includegraphics[width=1.1in,height=1.3in]
		{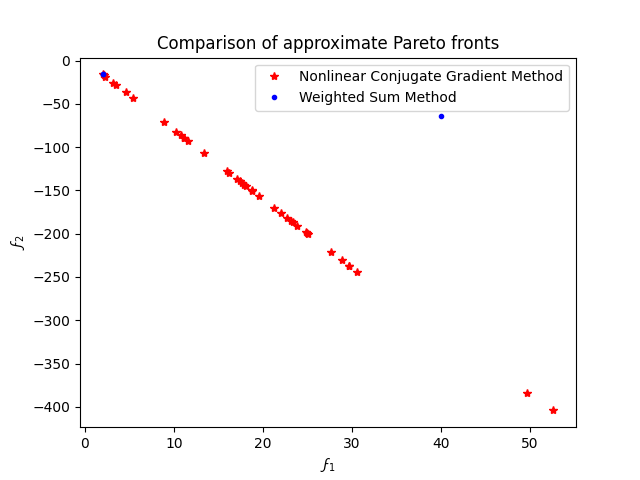}
		\centering  \caption*{($a_6$)}
	\end{subfigure}%
	\begin{subfigure}{1.2in}
		\includegraphics[width=1.1in,height=1.3in]
		{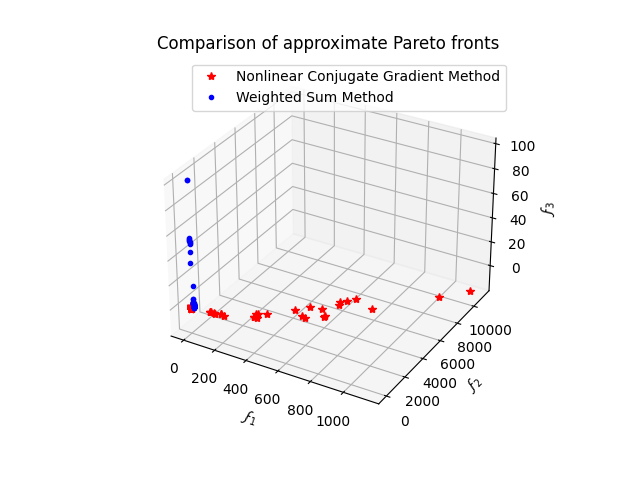}
		\centering    \caption*{($a_7$)}
	\end{subfigure}%
	\begin{subfigure}{1.2in}
		\includegraphics[width=1.1in,height=1.3in]
		{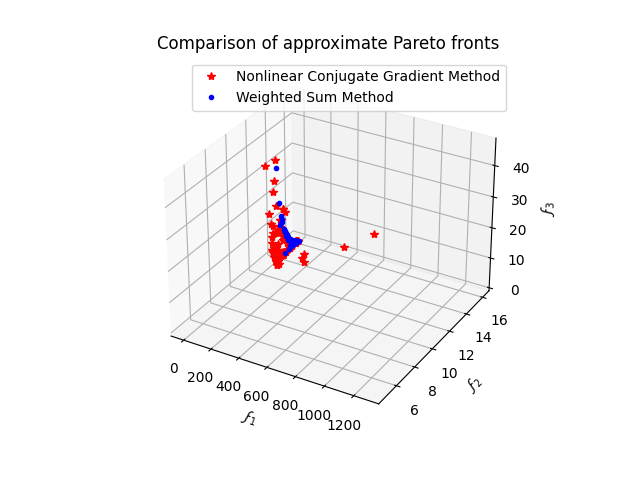}
		\centering     \caption*{($a_8$)}
	\end{subfigure}\\
	\begin{subfigure}{1.2in}
		\includegraphics[width=1.1in,height=1.3in]
		{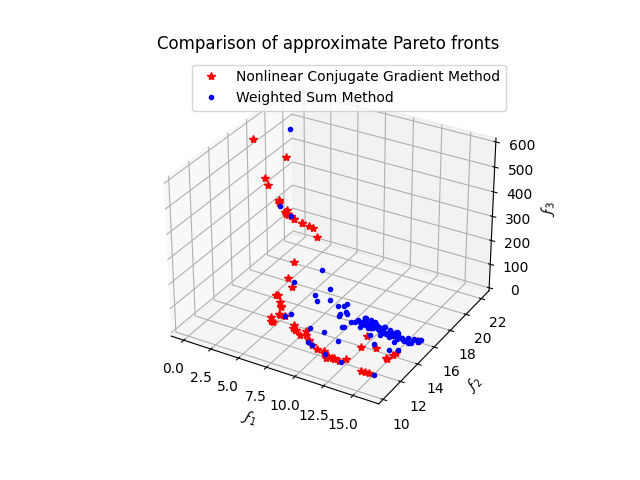}
		\centering    \caption*{($a_9$)}
	\end{subfigure}
	\begin{subfigure}{1.2in}
		\includegraphics[width=1.1in,height=1.3in]
		{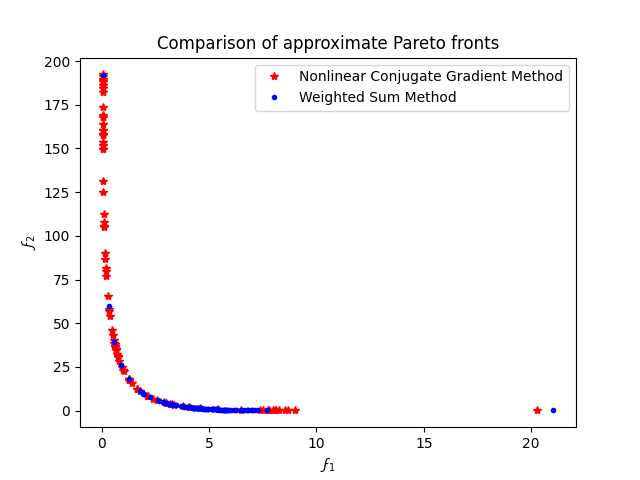}
		\centering   \caption*{($a_{10}$)}
	\end{subfigure}
	\begin{subfigure}{1.2in}
		\includegraphics[width=1.1in,height=1.3in]
		{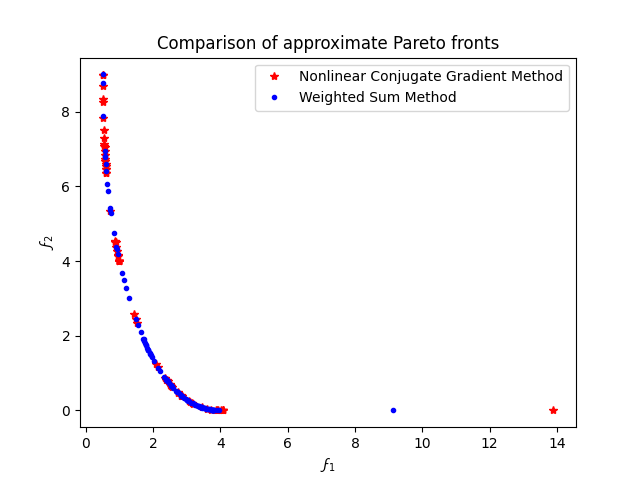}
		\centering  \caption*{($a_{11}$)}
	\end{subfigure}
	\begin{subfigure}{1.2in}
		\includegraphics[width=1.1in,height=1.3in]
		{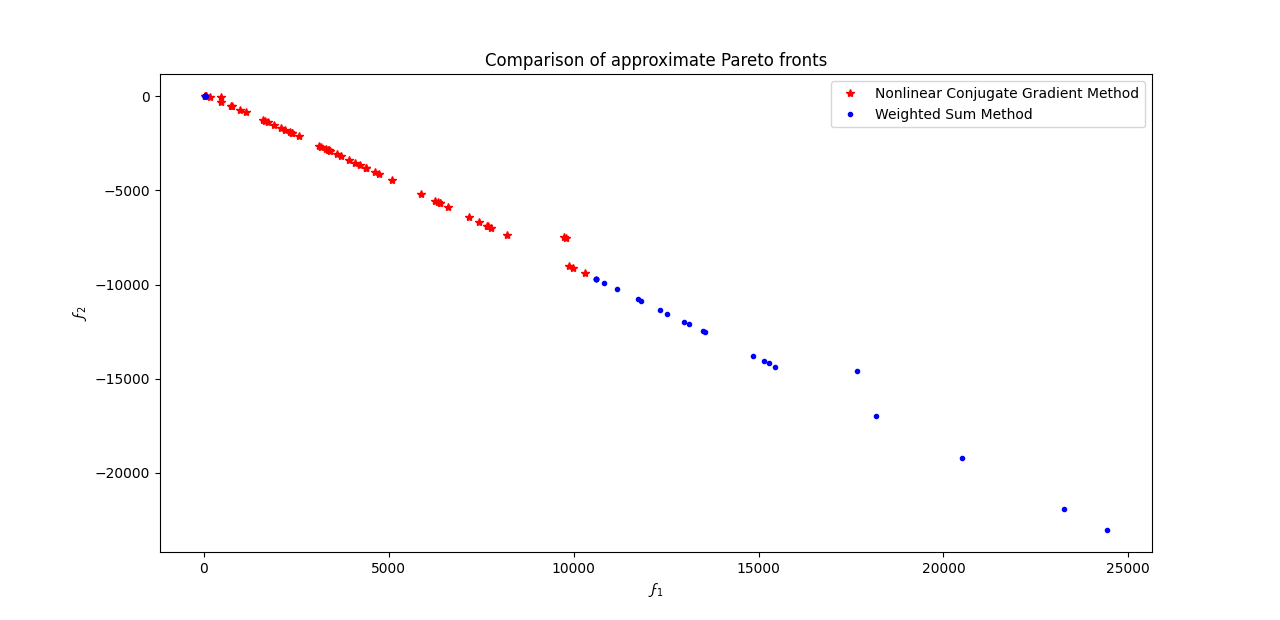}
		\centering  \caption*{($a_{12}$)}
	\end{subfigure}\\
	\begin{subfigure}{1.2in}
		\includegraphics[width=1.1in,height=1.3in]
		{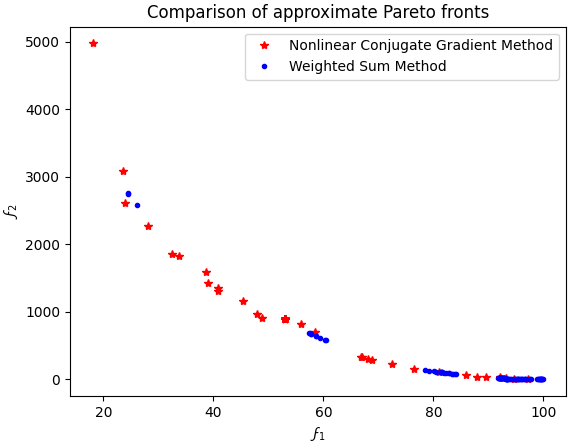}
		\centering  \caption*{($a_{13}$)}
	\end{subfigure}%
	\begin{subfigure}{1.2in}
		\includegraphics[width=1.1in,height=1.3in]
		{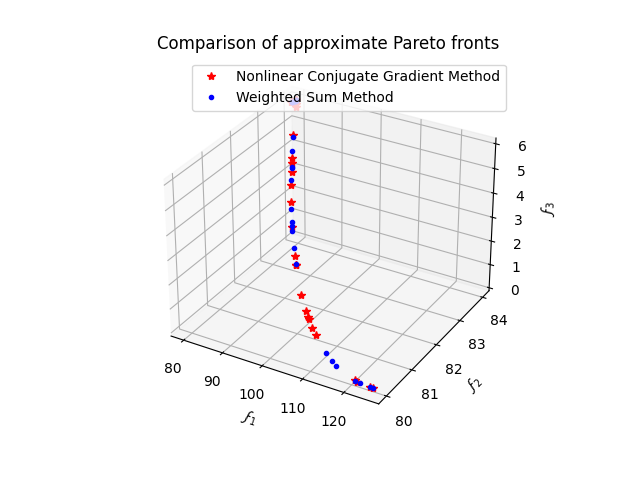}
		\centering  \caption*{($a_{14}$)}
	\end{subfigure}
	\begin{subfigure}{1.2in}
		\includegraphics[width=1.1in,height=1.3in]
		{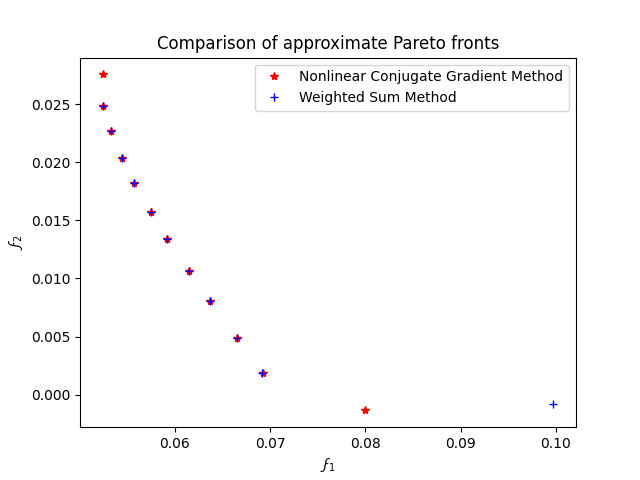}
		\centering  \caption*{($a_{15}$)}
	\end{subfigure}
	\begin{subfigure}{1.2in}
		\includegraphics[width=1.1in,height=1.3in]
		{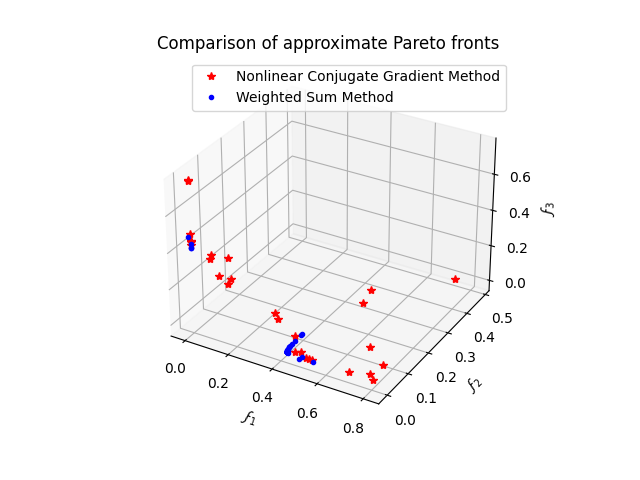}
		\centering  \caption*{($a_{16}$)}
	\end{subfigure}\\
	\begin{subfigure}{1.2in}
		\includegraphics[width=1.1in,height=1.3in]
		{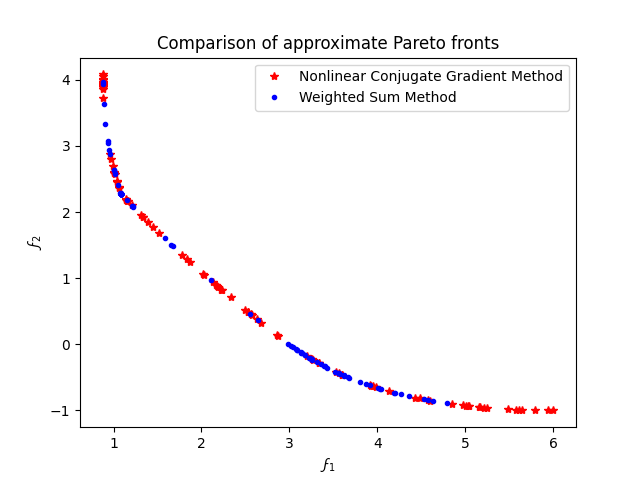}
		\centering  \caption*{($a_{17}$)}
	\end{subfigure}
	\begin{subfigure}{1.2in}
		\includegraphics[width=1.1in,height=1.3in]
		{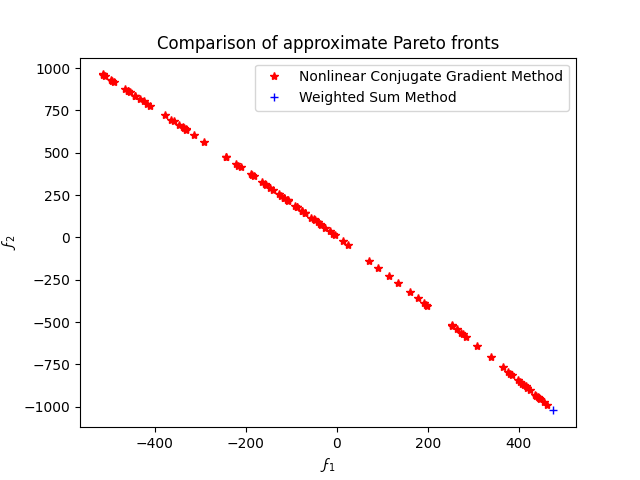}
		\centering  \caption*{($a_{18}$)}
	\end{subfigure}
	\begin{subfigure}{1.2in}
		\includegraphics[width=1.1in,height=1.3in]
		{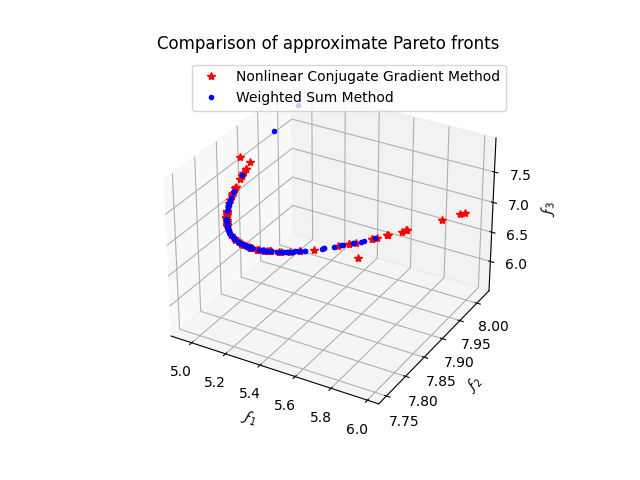}
		\centering  \caption*{($a_{19}$)}
	\end{subfigure}
	\begin{subfigure}{1.2in}
		\includegraphics[width=1.1in,height=1.3in]
		{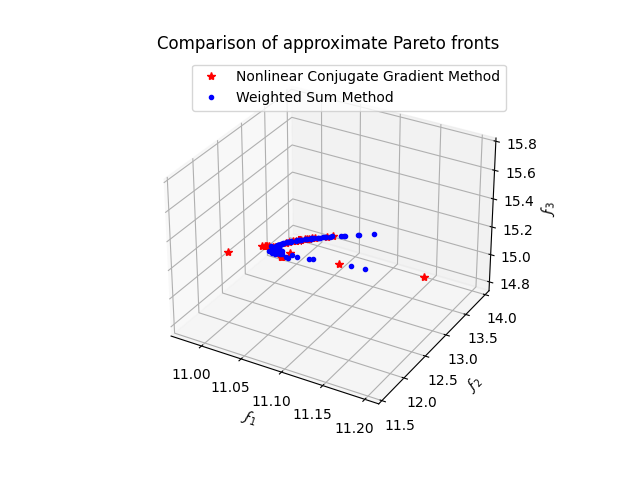}
		\centering  \caption*{($a_{20}$)}
	\end{subfigure}
	\caption{Comparison of approximate Pareto fronts generated by nonlinear conjugate gradient method (Algorithm \ref{algo1}) and weighted sum method for the test problems (TP1-TP20).}  
	\label{fig1}
\end{figure}

 We now examine 20 test problems listed in Table \ref{table1}, with details provided in Appendix 1. The full numerical solution for problem TP1 using Algorithm \ref{algo1} is presented in Example \ref{exmp1}, and solutions for the remaining problems can be obtained similarly.
\begin{table}[!ht]
	\caption{Details of test problems}
	\label{table1}
	\begin{center}
		\begin{tabular}
			{|c|c|c|c|}\hline
			Problem & $(m,n,p)$& $lb^T$ & $ub^T$\\ \hline
			
			
			TP1 &	$(2,1,2)$&$-5$&$5$ \\ \hline
			TP2 &$(2,2,2)$&$(-4,-4)^T$&$(4,4)^T$\\ \hline
			TP3 &$(2,3,3)$&$(0,0,0)^T$&$(1,1,1)^T$\\ \hline
			TP4 & (3,3,3)& $(1,-2,0)^T$ & $(3.5,2,1)^T$ \\ \hline
			TP5 & (2,2,2)& $(-6,-6)^T$ & $(6,4)^T$ \\ \hline
			TP6&$(2,1,2)$&$-3$&$3$\\ \hline
			TP7&$(3,3,3)$&$(-1,-1,-1)^T$&$(5,5,5)^T$\\ \hline
			TP8&$(3,2,3)$&$(-1,-1)^T$&$(5,2)^T$\\ \hline
			TP9&$(3,2,3)$&$(-1,-1)^T$&$(0,0)^T$\\ \hline
			TP10&(2,2,2)&$(-2,-2)^T$&$(5,5)^T$\\ \hline
			TP11&(2,2,2)&$(-6,-6)^T$&$(6,4)^T$\\ \hline
			TP12&(2,1,2)&$-100$&$100$\\ \hline
			TP13&(2,2,2)&$(0,0)^T$&$(1.0,1.0)^T$\\ \hline
			TP14&(3,2,3)&$(1,1)^T$&$(3,3)^T$\\ \hline
			TP15&(2,2,3)&$(.001,.001)^T$&$(1,1)^T$\\ \hline
			TP16&(3,10,3)&$(.001,...,.001)^T$&$(1,...,1)^T$\\ \hline
			TP17&(2,2,2)&$(-4,-4)^T$&$(5,5)^T$\\ \hline
			TP18&(2,1,2)&$-6$&$6$\\ \hline
			TP19&(3,3,3)&$(1,-2,0)^T$&$(3.5,2,1)^T$\\ \hline
			TP20&(3,3,3)&$(-1,-2,-1)^T$&$(4,5,3.4)^T$\\ \hline
		\end{tabular}
		
	\end{center}
\end{table}
\begin{example} (Nonconvex problem)\label{exmp1}
	Consider the problem
	$$\min_{x\in\mathbb{R}}\left(\digamma_1(x),\digamma_2(x)\right),$$
	where $\digamma_1(x)= \max\{h_{1}(x,\omega^{i}): i=1,2\}$, $\digamma_2(x)=\max\{h_{2}(x;\omega^{i}): i=1,2\}$, $\omega^1=-1, \omega^2=3.$ Note that $h_{1}(x,\omega^1)=(x+1)^2$, $h_{1}(x,\omega^2)=(x-3)^2$, $h_{2}(x,\omega^1)=x^2+x$ and $h_{2}(x,\omega^{2})=x^2-3x$. \\ \\
	We consider $lb=-5$ and $ub=5$ to generate an approximate Pareto front. Approximate Pareto fronts in both methods are given in Figure~\ref{fig1}~$(a_1)$.
	Consider $x^0= -0.6310622,$ $h_{1}(x^0,\omega^1)=0.1361151,~h_{1}(x^0,\omega^2)=13.18461274,$
	$h_{2}(x^0,\omega^1)=-1.02930171,~ h_{2}=(x^0,\omega^2)=1.49494711,$ $\digamma_1(x^0)=13.18461274,  ~\digamma_2(x^0)=1.49494711,$ 
	$v^0=s(x^0)= 0.57929217,$ $t=-3.36e-01,$ $\alpha_0=1,$
	$x^1=x^0+\alpha_0d^0=-0.05177003.$ At $x^1$, $h_{1}(x^1,\omega^1)=2.66036355,~h_{1}(x^1,\omega^2)=5.61186681,~ 
	h_{2}(x^1,\omega^1)=0.23282273,~h_{2}(x^1,\omega^2)= -2.29142564,~ \digamma_1(x^1)= 5.61186681,~\digamma_1(x^2)=0.23282273.$ We can observe that $\digamma_1(x^1)<\digamma_1(x^0),$ $\digamma_2(x^1)<\digamma_2(x^0).$ Therefore, $\digamma(x^1)<\digamma(x^0)$ and hence $v^0$ is a descent direction for $\digamma$ at $x^0.$
	Now, $s(x^1) = (0.63106215)$, $\beta_0=1$ $v^1=s(x^1)+\beta_0 v^0=0.63106215+0.57929217=1.21035432.$ Now we calculate $x^2=x^1+\alpha_1v^1=0.55340713,$ where $\alpha_1=0.5.$ Therefore, at $x^2,$  $h_{1}(x^2,\omega_1)=4.6453634,~  h_{1}(x^2,\omega_2)= 3.40287931,$
	$h_{2}(x^2,\omega_1)=-0.17943187,~ h_{2}(x^1,\omega_2)=-4.80067391,$ $\digamma_1(x^2)= 4.6453634,~ \digamma_1(x^2)=-0.17943187.$ We can observe that $\digamma_1(x^2)<\digamma_1(x^1),$ $\digamma_2(x^2)<\digamma_2(x^1),$ and we obtained $\digamma(x^2)<\digamma(x^2).$ Therefore, $v^1$ is a descent direction for $\digamma$ at $x^1.$ After 3 iteration, we get $T(x^*)=t= -3.16e-09$ and $s(x^*) =-1.38832447e-09,$ which satisfies the stopping criterion of the algorithm, and $x^*=1.15531051.$ Also, it is observed that $0\in conv \{\partial \digamma_{1}(x^*),\partial \digamma_{2}(x^*)\}.$ Therefore, by Lemma \ref{lm1}, $x^*$ is critical point for $\digamma$.

\end{example}

The figures, labeled from Figure \ref{fig1}($a_1$) to Figure \ref{fig1}($a_{20}$), illustrate the comparison between the approximate Pareto front generated by the NCGM and the weighted sum method for each of the test problems TP1-TP20.  
\par
\begin{figure}[H]
	\begin{subfigure}{3.0in}
		\includegraphics[width=2.6in,height=2.0in]
		{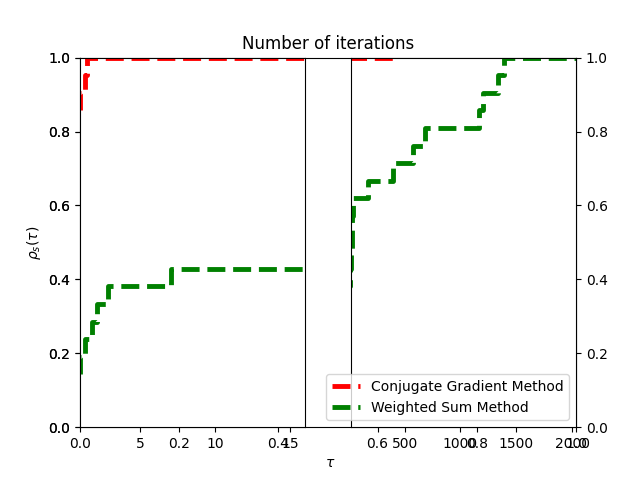}
		\centering    \caption*{($a_1$)}
	\end{subfigure}%
	\begin{subfigure}{3.0in}
		\includegraphics[width=2.6in,height=2.0in]
		{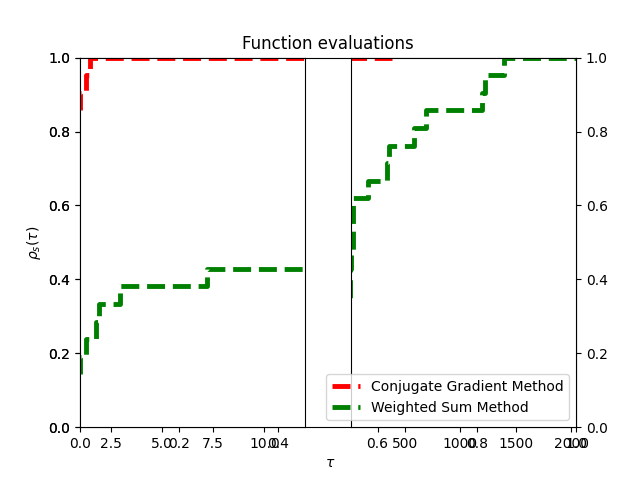}
		\centering     \caption*{($a_2$)}
	\end{subfigure}\\
	\begin{subfigure}{3.0in}
		\includegraphics[width=2.6in,height=2.0in]
		{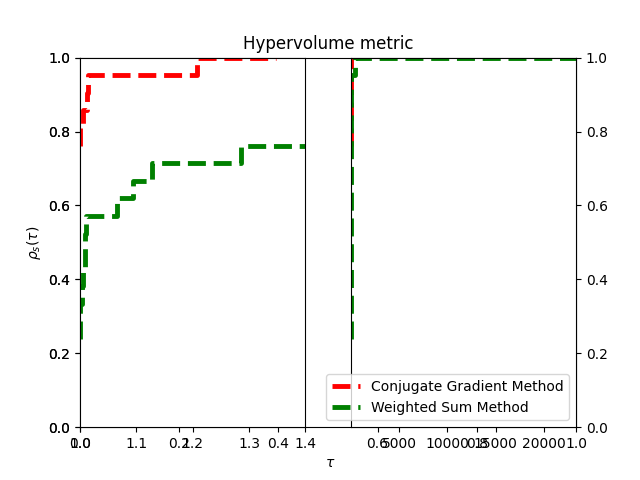}
		\centering   \caption*{($a_3$)}
	\end{subfigure}
	\begin{subfigure}{3.0in}
		\includegraphics[width=2.6in,height=2.0in]
		{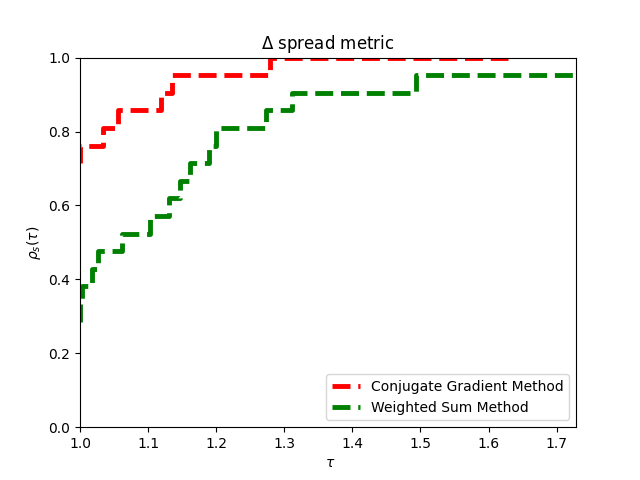}
		\centering   \caption*{($a_4$)}
	\end{subfigure}
	\caption{Comparison of the performance profiles of the number of iteration, function evaluation,  hypervolume metric, $\Delta$ spread metric for nonlinear conjugate gradient method and weighted sum method.}  
	\label{fig3}
\end{figure}

\textbf{Performance profiles:} The weighted sum method and Algorithm \ref{algo1} are compared using performance profiles.~Performance profiles are employed to contrast various approaches (see \citet{ansary2015modified,kumar2024steepest} for more details of performance profiles). The numerical results will be shown using performance profiles graphics, which are useful tools for comparing several methods on a large set of test problems. Let $\mathcal{SO}$ be the set of solvers or methods, $\mathcal{P}$ be the set of problems, and $\varsigma_{p,s} > 0$ be the performance of the solver $s \in \mathcal{SO}$ on the problem $p \in \mathcal{P}$, where lower values of $\varsigma_{p,s}$ mean better performances. Define the performance ratio $r_{p,s} := \frac{\varsigma_{p,s}}{\min\{\varsigma_{p,s} | s \in \mathcal{SO}\}}$. Then, the performance profile is obtained by plotting, for all $s \in\mathcal{SO}$, the cumulative distribution function $\rho_s : [1,\infty[ \rightarrow [0, 1]$ for the performance ratio $r_{p,s}$ given by $\rho_s(\tau) := \frac{1}{|\mathcal{P}|}|\{p \in \mathcal{P} | r_{p,s} \leq \tau \}|$, where $|\cdot|$ denotes the cardinality of the set. In a performance profile graphic, $\rho_s(\tau = 1)$ is the fraction of problems for which solver $s$ was the most efficient over all the methods. On the other hand, $\rho_s(\tau \equiv \infty)$ represents the fraction of problems for which solver $s$ was able to find a solution, independently of the required effort. Therefore, the fractions $\rho_s(\tau = 1)$ and $\rho_s(\tau \equiv \infty)$, which can be accessed on the extreme left and right of the graph, are usually associated with the efficiency and robustness of solver $s$, respectively.\\
$~~~~$To justify how much well-distributed this set is, the following metrics are considered for computing performance profile.
In multiobjective optimization, we are mainly interested in estimating the Pareto frontier of a given problem. A commonly used strategy for this task is to run an algorithm from several starting points and collect the efficient points found. We compare the results using the well known $\Delta$-spread and Hypervolume metrics.\\
\textbf{$\Delta$-spread metric:} $\Delta$-spread metric measures the ability to obtain well-distributed points along the Pareto frontier. Let $x^1,x^2,\dots,x^N$ be the set of points obtained by a solver $s$ for problem $p$ and let these points be sorted by $f_j(x^i)\leq f_j(x^{i+1})$ $(i=1,2,\dots,N-1)$. Suppose $x^0$ is the best known approximation of global minimum of $f_j$ and $x^{N+1}$ is the best known global maximum of $f_j$, computed over all the approximated Pareto fronts obtained by different solvers. Define $\bar{\delta_j}$ as the average of the distances $\delta_{i,j}$, $ i=1,2,\dots,N-1.$ For an algorithm  $s$ and a problem $p$, the spread metric $\Delta_{p,s}$ is
\begin{eqnarray*}
	~~~~~~~~~~~~~~~~~~~~~~	\Delta_{p,s}:=\underset{j\in\Lambda_m}{\max}\left(\frac{\delta_{0,j}+\delta_{N,j}+\Sigma_{i=1}^{N-1} |\delta_{i,j}-\bar{\delta_j}|}{\delta_{0,j}+\delta_{N,j}+(N-1)\bar{\delta_j}}\right).
\end{eqnarray*}
\textbf{Hypervolume metric:} Hypervolume metric of an approximate Pareto front with respect to a reference point $P_{ref}$ is defined as the volume of the total region dominated by the efficient solutions obtained by a method with respect to the reference point. We have used the codes from \url{https://github.com/anyoptimization/pymoo} to calculate hypervolume metric. Higher values of $hv_{p,s}$ indicate better performance using hypervolume metric. So while using the performance profile of the solvers measured by hypervolume metric we need to set $\widetilde{hv_{p,s}}=\frac{1}{hv_{p,s}}$.\\
Performance profile using hypervolume metric and $\Delta$ spread metric are given in Figure \ref{fig3}($a_3$) and Figure \ref{fig3}($a_4$), respectively. 
Regarding hypervolume metric (resp. $\Delta$ spread metric), the
efficiencies of the algorithms are $76.0\%$ and $24.0\%$ (resp. $71.2\%$ and $28.11\%$) for Algorithm \ref{algo1} and the weighted sum method, respectively.
\par In addition to this, we have computed performance profiles using number of iterations and function evaluations. Gradients in both methods are calculated using forward difference formula which requires $n$ additional function evaluations. If {\it \# Iter} and {\it \# Fun} denote the number of iterations and number of function evaluations required to solve a problem respectively then total function evaluation is ${\it \# Fun}+n  {\it \# Iter}$. Performance profiles using number of iterations and number of function evaluations are given in Figure \ref{fig3}($a_1$) and Figure \ref{fig3}($a_2$), respectively. Regarding function evaluation (resp. number of iterations), the efficiencies of the algorithms are $85.0\%$ and $14.9\%$ (resp. $86.6\%$ and $15.7\%$) for Algorithm \ref{algo1} and the weighted sum method, respectively. One can observe from performance profile figures that Algorithm \ref{algo1} performs better than weighted sum method in most cases.
\section{Conclusion}\label{sec6}
We solved an uncertain multiobjective optimization problem $P(U^F)$ defined in (\ref{main uncertain mop}), using its robust counterpart $OWC_{P(U^F)}$, which is a deterministic multiobjective optimization problem. By solving $OWC_{P(U^F)}$, we obtain the solution for $P(U^F)$, eliminating the need to directly solve $P(U^F)$. To tackle $OWC_{P(U^F)}$, we employed a nonlinear conjugate gradient method. Our approach involved developing a nonlinear conjugate gradient algorithm to find a critical point. This algorithm involves a subproblem to determine the steepest descent direction at the current iteration. By using the descent direction at the current iteration and the previous direction a new direction is updated. To update the new direction, Fletcher-Reeves, conjugate descent, Dai-Yuan, Polak-Ribi$\grave{e}$re-Polyak, and Hestenes-Stiefel parameters extension from deterministic MOP to UMOP were considered. For these parameters, it is shown that the updated direction satisfies the sufficient descent property (i.e., the updated direction is a descent direction).~To find the step length, a method for developing inexact line searches similar to Armijo's is created. A nonlinear conjugate gradient algorithm is written with the help of an updated direction and step length size. The convergence analysis of the proposed method is also discussed. At the end of the paper, some test problems are constructed to validate the nonlinear conjugate gradient algorithm. A comparison of the nonlinear conjugate gradient method with the weighted sum method is also presented by using performance profiles. We also compared the approximate Pareto front obtained by the nonlinear conjugate gradient method with the one generated by the weighted sum method. We observed that the nonlinear conjugate gradient method successfully generates good approximations of the Pareto fronts in both convex and nonconvex cases.~However, the weighted sum method fails to produce even an approximate front in the nonconvex case. In addition to this, we have computed performance profiles using several iterations, function evaluations, $\Delta$ spread metric, and hypervolume metric. Based on performance profile figures the nonlinear conjugate gradient method better than the weighted sum method in most cases. The choice of weights, constraints, or the importance of the functions, respectively, is not known in advance and needs to be pre-specified, which is an inherent drawback of the scalarization method for uncertain multiobjective optimization problems. In addition, the nonlinear conjugate gradient method for $OWC_{P(U^F)}$ does not require predetermined weighting factors or any other kind of predetermined ranking or ordering information for objective functions, eliminating another drawback of scalarization methods.
\par For future work, we will develop the spectral conjugate gradient method, which requires less computational work and outperforms sophisticated conjugated methods in many problems.
\section*{Availability of data and materials}
Not applicable.
\section*{Acknowledgment:}
This research is supported by Govt. of India CSIR fellowship, Program No. 09/1174(0006)/2019-EMR-I, New Delhi.
\section*{Conflict of interest:}
The authors say they have no competing interests.
\bibliographystyle{apa}
\bibliography{Conjugate_references}

\newpage
\newpage
\appendix
\section{Details of test problems}
\begin{appendices}{\bf Appendix A: Details of test problems}
	\begin{enumerate}[({TP}1)]
		\item $\digamma: \mathbb{R}\rightarrow \mathbb{R}^2$ is defined as $\digamma_j(x)= \max\{h_{j}(x,\omega^{i}): i=1,2\}$ $j=1,2$ where $\omega^1=-1$, $\omega^2=3$, $h_{1}(x,\omega^i)=(x-\omega^i)^2$ and $h_{2}(x,\omega^i)=x^2+\omega^i x$ for $i=1,2$.
		\item $\digamma: \mathbb{R}^2\rightarrow \mathbb{R}^2$ is defined as $\digamma_j(x)= \max\{h_{j}(x,\omega^{i}): i=1,2\}$ $j=1,2$ where $\omega^1=(1,3)^T,~ ~\omega^2=(3,1)^T.$
		Note that
		$h_{1}(x,\omega^1)=(x_1-1)^2+(x_2-3)^2$, and $h_{1}(x,\omega^2)=(x_1-3)^2+(x_2-1)^2$, $h_{2}(x,\omega^1)=x_1^2+3x_2^2$, and $h_{2}(x,\omega^{2})=3x_1^2+x_2^2$.
		\item $\digamma: \mathbb{R}^3\rightarrow \mathbb{R}^2$ is defined as $\digamma_j(x)= \max\{h_{j}(x,\omega^{i}): i=1,2,3\},$ $j=1,2.$ Here $\omega^i=(\omega^i_1,\omega^i_2,\omega^i_3)$ and therefore $\omega^1=(1,1,1)^T,$  $\omega^2=(1,-1,1)^T$, $\omega^3=(1,-2,2)^T$,
		\begin{eqnarray*}
			h_{1}(x,\omega^i)=&1-e^{-\sum_{j=1}^3\omega^j_{i} \left(x_j-\frac{1}{\sqrt{3}}\right)^2}~~i=1,2,3\\ \mbox{and}~~~~
			h_{2}(x,\omega^i)=&1-e^{-\sum_{j=1}^3\omega^j_{i} \left(x_j+\frac{1}{\sqrt{3}}\right)^2}~~i=1,2,3.
		\end{eqnarray*}
		\item $\digamma: \mathbb{R}^3\rightarrow \mathbb{R}^3$ is defined as $\digamma_j(x)= \max\{h_{j}(x,\omega^{i}): i=1,2,3\}$ $j=1,2,3,$ where $\omega^i=(\omega^i_1,\omega^i_2,\omega^i_3)$ and
		$\omega^1=(1,1,1)^T$, $\omega^2=(1,-1,1)^T$, $\omega^3=(1,-2,2)^T.$
		\begin{eqnarray*}
			h_{1}(x,\omega^i)&=&\left(1+\omega^i_3x_3\right)\left(\omega^i_{1}\omega^i_{2} x_1^3x_2^3-10\omega^i_{1}x_1-4\omega^i_{2}x_2\right)~~i=1,2,3\\
			h_{2}(x,\omega^i)&=&\left(1+\omega^i_3x_3\right)\left(\omega^i_{1}\omega^i_{2}x_1^3x_2^3-10\omega^i_{1}x_1+4\omega^i_{2}x_2\right)~~i=1,2,3\\
			\mbox{and}~~~~
			h_{3}(x,\omega^i)&=&\left(1+\omega^i_{3}x_3\right)\omega^i_{1}x_1^2~~i=1,2,3.
		\end{eqnarray*}
		
		\item\label{T5}  $\digamma: \mathbb{R}^2\rightarrow \mathbb{R}^2$ is defined as $\digamma_j(x)= \max\{h_{j}(x,\omega^{i}): i=1,2\},$ $j=1,2,$ where $\omega^i=(\omega^i_1,\omega^i_2)$ and $\omega^1=\{(2,2)$, $\omega^2=(0,4).$
		\begin{eqnarray*}
			h_{1}(x,\omega^i)&=&(x_1-\omega^i_1)^2+(x_2+\omega^i_2)^2~~i=1,2\\ \mbox{and}~~~~
			h_{2}(x,\omega^i)&=&(\omega^i_1x_1+\omega^i_2x_2)^2~~i=1,2.
		\end{eqnarray*}
		\item $\digamma: \mathbb{R}\rightarrow \mathbb{R}^2$ is defined as $\digamma_j(x)= \max\{h_{j}(x,\omega^{i}): i=1,2\},$ $j=1,2,$ where $\omega^1=-2$, $\omega^2=5.$
		\begin{eqnarray*}
			h_{1}(x,\omega^i)&=&(x-\omega^i)^2~~i=1,2\\ \mbox{and}~~~~
			h_{2}(x,\omega^i)&=&-x^2-\omega^i x~~i=1,2.
		\end{eqnarray*}
		\item $\digamma: \mathbb{R}^3\rightarrow \mathbb{R}^3$ is defined as $\digamma_j(x)= \max\{h_{j}(x,\omega^{i}): i=1,2,3\},$ $j=1,2,3,$ where $\omega^i=(\omega^i_1,\omega^i_2)$ and $\omega^1=(4,1),$ $\omega^2=(0,2)$, and $\omega^3=(1,0).$	
		\begin{eqnarray*}
			h_{1}(x,\omega^i)&=&x_1^2+(x_2-\omega^i_1)^2-\omega^i_2x_3^2~~i=1,2,3\\
			h_{2}(x,\omega^i)&=&\omega^i_1x_1+\omega^i_2x_2^2+x_3+4\omega^i_1\omega^i_2~~i=1,2,3\\\mbox{and}~~~~
			h_{3}(x,\omega^i)&=&\omega^i_1x_1^2+6x_2^2+25(x_3-\omega^i_2x_1)^2.
		\end{eqnarray*}
		\item $\digamma: \mathbb{R}^2\rightarrow \mathbb{R}^3$ is defined as $\digamma_j(x)= \max\{h_{j}(x,\omega^{i}): i=1,2,3\},$ $j=1,2,3,$ where $\omega^i=(\omega^1_i,\omega^2_i)$ and $\omega^1=(2,3),$ $\omega^2=(4,5)$, and $\omega^3=(2,0).$
		\begin{eqnarray*}
			h_{1}(x,\omega^i)&=& x_1^2+\omega^i_1x_2^4+\omega^i_1\omega^i_2x_1x_2,~~i=1,2,3\\
			h_{2}(x,\omega^i)&=&5x_1^2+\omega^i_1x_2^2+\omega^i_2x_1^4x_2,~~i=1,2,3\\\mbox{and}~~~~
			h_{3}(x,\omega^i)&=&e^{-\omega^i_1x_1+\omega^i_2x_2}+x_1^2-\omega^i_1x_2^2,~~i=1,2,3.
		\end{eqnarray*}
		\item  $\digamma: \mathbb{R}^2\rightarrow \mathbb{R}^3$ is defined as $\digamma_j(x)= \max\{h_{j}(x,\omega^{i}): i=1,2,3\},$ $j=1,2,3,$ where $\omega^i=(\omega^i_i,\omega^i_2)$ and $\omega^1=(2,3),$ $\omega^2=(1,2)$, and $\omega^3=(4,5).$
		\begin{eqnarray*}
			h_{1}(x,\omega^i)&=&100\omega^i_1 (x_2-x_1^2)^2+\omega^i_2(1-x_1)^2~~i=1,2,3\\
			h_{2}(x,\omega^i)&=&(x_2-\omega^i_1)^2+\omega^i_2x_1^2~~i=1,2,3\\\mbox{and}~~~~
			h_{3}(x,\omega^i)&=&\omega^i_1x_1^2+3\omega^i_2x_2^2 ~~i=1,2,3.
		\end{eqnarray*}
		\item\label{T10} $\digamma: \mathbb{R}^2\rightarrow \mathbb{R}^2$ is defined as $\digamma_j(x)= \max\{h_{j}(x,\omega^{i}): i=1,2\},$ $j=1,2,$ where $\omega^i=(\omega^i_1,\omega^i_2,\omega^i_3)$ and $\omega^1=(1,2,2),$ $\omega^2=(1,3,0).$
		\begin{eqnarray*}
			h_{1}(x,\omega^i)&=&\omega^i_1x_1^2+\omega^i_2x_2^2+\omega^i_1x_1+\omega^i_1\omega^i_3x_2,~~i=1,2\\\mbox{and}~~~~h_{2}(x,\omega^i)&=&(\omega^i_1+\omega^i_2x_2)^2+\omega^i_1x_1+x_2+10(x_1+\omega^i_3x_2)+e^{(1+\omega^i_1x_1+\omega^i_2x_2)^2}~~i=1,2.
		\end{eqnarray*}
		{\bf Note:} TP\ref{T5}-TP\ref{T10} are taken from the \citet{shubham2023newton}.
		\item  $\digamma: \mathbb{R}^2\rightarrow \mathbb{R}^2$ is defined as $\digamma_j(x)= \max\{h_{j}(x,\omega^{i}): i=1,2\},$ $j=1,2,$ where $\omega^i=(\omega^1_i,\omega^2_i)$ and $\omega^1=(1,2),$ $\omega^2=(2,3).$
		\begin{eqnarray*}
			\text{and}~~	h_{1}(x,\omega^i)&=&\frac{1}{4}(x_1-\omega^i_1)^4+2(x_2-\omega^i_2)^4~~i=1,2\\\mbox{and}~~~~	h_{2}(x,\omega^i)&=&(\omega^i_1x_2-\omega^i_2x_1^2)^2+(1-\omega^i_1x_1)^2,~i=1,~2.
		\end{eqnarray*}
		\item 	$\digamma: \mathbb{R}\rightarrow \mathbb{R}^2$ is defined as $\digamma_j(x)= \max\{h_{j}(x,\omega^{i}): i=1,2\},$ $j=1,2,$ where $\omega^1=-3,$ $\omega^2=8.$
		\begin{eqnarray*}
			h_1(x,\omega^i)&=&(x-\omega^i)^2,~~i=1,2\\\mbox{and}~~
			h_{2}(x,\omega^i)&=&-x^2-\omega^i x.
		\end{eqnarray*}
		\item $\digamma: \mathbb{R}^2\rightarrow \mathbb{R}^2$ is defined as $\digamma_j(x)= \max\{h_{j}(x,\omega^{i}): i=1,2\},$ $j=1,2,$ where $\omega^i=(\omega^i_1,\omega^i_2)$ and $\omega^1=(1,1),$ $\omega^2=(0,2).$
		\begin{eqnarray*}
			h_{1}(x,\omega^i)&=&(x_1-\omega^i_1)^2+(x_2+\omega^i_2)^2,~~i=1,~2\\\mbox{and}~~
			h_{2}(x,\omega^i)&=&(\omega^i_1x_1+\omega^i_2x_2)^2~~i=1,~2.
		\end{eqnarray*}
		\item $\digamma: \mathbb{R}^2\rightarrow \mathbb{R}^3$ is defined as $\digamma_j(x)= \max\{h_{j}(x,\omega^{i}): i=1,2,3\},$ $j=1,2,3$ where $\omega^i=(\omega^i_1,\omega^i_2)$ and $\omega^1=(4,1),$ $\omega^2=(5,2),$ and $\omega^3=(6,4).$
		\begin{eqnarray*}
			h_{1}(x,\omega^i)&=&100\omega^i_1 (x_2-x_1^2)^2+\omega^i_2(1-x_1)^2,~~i=1,2,3\\
			h_{2}(x,\omega^i)&=&(x_2-\omega^i_1)^2+\omega^i_2x_1^2~~i=1,2,3\\ \mbox{and}~~
			h_{3}(x,\omega^i)&=&\omega^i_1x_1^2+3\omega^i_2x_2^2~~i=1,2,3.
		\end{eqnarray*}
		
		\item \label{dt1} $\digamma: \mathbb{R}^n\rightarrow \mathbb{R}^m$ is defined by $\digamma_j(x)= \max\{h_{j}(x,\omega^{i}): i=1,2,3\},$ $j=1,2,\dots,m,$ where $\omega^1=(0.25,0.25,\dots,0.25)^T\in \mathbb{R}^n$,  $\omega^2=(0.5,0.5,\dots,0.5)^T\in\mathbb{R}^n$, $\omega^3=(0.75,0.75,\dots,0.75)^T\in\mathbb{R}^n$. Define $$gx_i=100\left(K+\sum_{k=m}^n\left((x_i-\omega^i_k)^2-\cos(20\pi(x_k-\omega^i_k))\right)\right)$$
		for $i=1,2,3$ where $K=m+n-1$. $f_j(x,\omega^i) $ for $j=1,2,..., m$ are defined as
		\begin{eqnarray*}
			h_{1}(x,\omega^i)&=&0.5(1+gx_i)\prod_{k=1}^{m-1} x_k\\
			h_j(x,\omega^i)&=&	0.5(1+gx_i)\prod_{k=1}^{m-j} x_k(1-x_{m-j+1}) ~~j=2,\dots,m
		\end{eqnarray*}
		\item \label{dt2} $\digamma: \mathbb{R}^n\rightarrow \mathbb{R}^m$ is defined by $\digamma_j(x)= \max\{h_{j}(x,\omega^{i}): i=1,2,3\},$ $j=1,2,\dots,m$ where $\omega^1=(0.4,0.4,\dots,0.4)^T\in \mathbb{R}^n$,  $\omega^2=(0.5,0.5,\dots,0.5)^T\in \mathbb{R}^n$, $\omega^3=(0.6,0.6,\dots,0.6)^T\in \mathbb{R}^n$. Define $$gx_i=\sum_{k=m}^n(x_i-\omega^i_k)^2$$
		for $i=1,2,3$. $f_j(x,\omega^i) $ for $j=1,2,..., m$ are defined as
		\begin{eqnarray*}
			h_{1}(x,\omega^i)&=&(1+gx_i)\prod_{k=1}^{m-1} \cos(0.5\pi x_k)\\
			f_j(x,\omega^i)&=&	0.5(1+gx_i)\prod_{k=1}^{m-j}\cos(0.5\pi x_k)\sin(0.5\pi x_{m-k+1}) ~~j=2,\dots,m
		\end{eqnarray*}
		\item $\digamma: \mathbb{R}^2\rightarrow \mathbb{R}^2$ is defined as $\digamma_j(x)= \max\{h_{j}(x,\omega^{i}): i=1,2\},$ $j=1,2,$ where $\omega^i=(\omega^i_1,\omega^i_2)$ and $\omega^1=(50,4),$ $\omega^2=(101,3).$
		\begin{eqnarray*}
			h_{1}(x,\omega^i)&=&\omega^i_1x_1^2+(\omega^i_2+\omega^i_1)x_2+\omega^i_1\omega^i_2x_1x_2+3~~i=1,~2\\\mbox{and}~~
			h_{2}(x,\omega^i)&=&\omega^i_1x_1^2+\omega^i_2x_2^2+(x_1+x_2)\omega^i_1\omega^i_2+4~~i=1,~2.
		\end{eqnarray*}
		\item $\digamma: \mathbb{R}\rightarrow \mathbb{R}^2$ is defined as $\digamma_j(x)= \max\{h_{j}(x,\omega^{i}): i=1,2\}$ $j=1,2$ where $\omega^1=-9$, $\omega^2=58$, $h_{1}(x,\omega^i)=-\omega^ix^2+57x+1$ and $h_{2}(x,\omega^i)=-\omega^ix^2-25x+4$ for $i=1,2$.
		\item $\digamma: \mathbb{R}^3\rightarrow \mathbb{R}^3$ is defined as $\digamma_j(x)= \max\{h_{j}(x,\omega^{i}): i=1,2,3\},$ $j=1,2,3$ where $\omega^i=(\omega^i_1,\omega^i_2,\omega^i_3)$ and $\omega^1=(76,4,4),$ $\omega^2=(0,9,6),$ and $\omega^3=(4,6,1).$
		\begin{eqnarray*}
			h_{1}(x,\omega^i)&=&(\omega^i_1+\omega^i_2)x_1^2+x_1^4x_2(\omega^i_1+\omega^i_3)+x_3\omega^i_1\omega^i_2+1~~~i=1,2,3\\
			h_{2}(x,\omega^i)&=&\omega^i_1\omega^i_3(x^2_1+2x^2_2)+(\omega^i_1+\omega^i_3)x^2_1x_2+3~~~i=1,2,3\\\mbox{and}~	h_{3}(x,\omega^i)&=&\omega^i_1x_1^3+\omega^i_2x_2^2+(\omega^i_1+\omega^i_2+\omega^i_3)x_1x_2x_3+6~~~i=1,2,3.\\
		\end{eqnarray*}
		\item $\digamma: \mathbb{R}^3\rightarrow \mathbb{R}^3$ is defined as $\digamma_j(x)= \max\{h_{j}(x,\omega^{i}): i=1,2,3\},$ $j=1,2,3$ where $\omega^i=(\omega^i_1,\omega^i_2,\omega^i_3)$ and $\omega^1=(1,0,90),$ $\omega^2=(9,17,6),$ and $\omega^3=(8,2,1).$
		\begin{eqnarray*}
			h_{1}(x,\omega^i)=\omega^i_1\omega^i_2\omega^i_3x_1^2+x_3^4(\omega^i_1+\omega^i_2)+x_1x_2(\omega^i_1\omega^i_2\omega^i_3)+\omega^i_1\omega^i_2\omega^i_3~~~i=1,2,3\\
			h_{2}(x,\omega^i)=x_1^3+(x_2^2+x_3^2)\omega^i_1\omega^i_2+\omega^i_1\omega^i_2(x_1^2+x_2^2+x_3^2)+\omega^i_1\omega^i_2~~~i=1,2,3\\\mbox{and}~	h_{3}(x,\omega^i)=x_1^2+x_2^2+(\omega^i_1+\omega^i_2+\omega^i_3)x_1x_2+\omega^i_1+1~~~i=1,2,3.\\
		\end{eqnarray*}
	\end{enumerate}
\end{appendices}
\end{document}